\documentclass[10pt,reqno]{amsart}
\usepackage[colorlinks=true,
pdfstartview=FitV, linkcolor=, citecolor=,
urlcolor=]{hyperref}
\usepackage{amsmath,amsfonts,latexsym,amssymb}
\usepackage{amssymb,amsfonts, amsmath}
\usepackage{mathrsfs}
\usepackage[latin1]{inputenc}
\usepackage[T1]{fontenc}
\usepackage{ae,aecompl}
\usepackage{braket}
\usepackage{comment}
\usepackage{amssymb,amsfonts, amsmath}
 \usepackage{mathtools} 
\usepackage[all,arc]{xy}
\usepackage{enumerate}
\usepackage{mathrsfs}
\usepackage{mathabx}
\usepackage{color}
\usepackage{graphicx}
\usepackage{subfig}
\usepackage{verbatim} 
\usepackage{amssymb,amsfonts, amsmath}
\usepackage[all,arc]{xy}
\usepackage{enumerate}
\usepackage{mathrsfs}
\usepackage{hyperref}
\usepackage{xstring}
\usepackage{thmtools}
\usepackage{thm-restate}

\newtheorem{theorem}{Theorem}[section]
\newtheorem{lemma}[theorem]{Lemma}
\newtheorem{proposition}[theorem]{Proposition}
\newtheorem{corollary}[theorem]{Corollary}
\newtheorem{definition}[theorem]{Definition}

\newtheorem{remark}[theorem]{Remark}
\newtheorem{ques}[theorem]{Question}
\newtheorem{conj}[theorem]{Conjecture}
\newtheorem{fact}[theorem]{Fact}

\newtheorem*{question}{Question}
\renewcommand{\leq}{\leqslant}
\renewcommand{\geq}{\geqslant}

\newcommand{\equaldef}{\overset{\mathrm{def}}{=}}

\long\def\@savemarbox#1#2{\global\setbox#1\vtop{\hsize\marginparwidth 
  \@parboxrestore\tiny\raggedright #2}}
\newcommand{\CC}{\mathbb C}

\newcommand{\Real}{\mathbb R}
\newcommand{\Z}{\mathbb Z}

\newcommand{\Sp}{\mathsf{Sp}}
\newcommand{\SU}{\mathsf{SU}}
\newcommand{\Isom}{\mathsf{Isom}}
\newcommand{\SL}{\mathsf{SL}}
\newcommand{\GL}{\mathsf{GL}}

\author[Tholozan]{Nicolas Tholozan}
\author[Tsouvalas]{Konstantinos Tsouvalas}
\title{Linearity and indiscreteness of amalgamated products of hyperbolic groups}

\date{\today}

\AtEndDocument{\bigskip{\footnotesize%
  \textsc{DMA-UMR8553, \'Ecole Normale Sup\'erieure, CNRS-PSL, Research University, 45, rue d'Ulm
75230 Paris Cedex 5
France} \par  
 \textit{E-mail address}: \texttt{nicolas.tholozan@ens.fr} 
\vspace{0.2cm}

\textsc{CNRS, Laboratoire Alexander Grothendieck, Institut des Hautes \'Etudes Scientifiques, 
Universite Paris-Saclay, 35 route de Chartres, 91440 Bures-sur-Yvette, France} \par  
 \textit{E-mail address}: \texttt{tsouvkon@ihes.fr}

  \textsc{}}}

\begin{document}

\begin{abstract} We discuss the linearity and discreteness of amalgamated products of linear word hyperbolic groups. In particular, we prove that the double of a torsion free Anosov group along a maximal cyclic subgroup is always linear, and we construct examples of such groups which do not admit any discrete and faithful representation into any simple Lie group of real rank $1$. We also build new examples of non-linear word hyperbolic groups, elaborating on a previous work of Canary--Stover--Tsouvalas.\end{abstract}

\frenchspacing
\maketitle
\section{Introduction}

In his groundbreaking work \cite{Gromov}, Gromov introduced the notion of \emph{word-hyperbolic group}, which captures the coarse geometric and algebraic properties of fundamental groups of closed negatively curved manifolds. Among their many interesting properties, hyperbolic groups are stable under many operations of topological origin, such as amalgamated products and HNN extensions over nice subgroups, see \cite{BF}. Moreover, certain quotients of a non-elementary hyperbolic group remain hyperbolic after adding sufficiently complicated relations (see \cite{Gromov}, \cite{Delzant} and \cite{Ol}).

The present work revolves around the following general question:

\begin{question}
Which word-hyperbolic groups admit geometric realisations ?
\end{question}

Of course, the term ``geometric realisation'' can be understood in many ways. Here, we will mainly be interested in realizing word-hyperbolic groups as discrete subgroups of (real) linear Lie groups.

An important family of word-hyperbolic groups of geometric origin is formed by convex-cocompact subgroups of simple rank one Lie groups. This includes linear groups that predate Gromov's definition by a century, such as  Schottky groups, quasi-Fuchsian groups or uniform hyperbolic lattices.\footnote{Gromov's theory was actually meant to generalize this class of groups in a coarse geometric context.} A convex-cocompact subgroup $\Gamma$ of a rank $1$ Lie group $G$ is quasi-isometrically embedded in the symmetric space $X$ of $G$, which is negatively curved, and the \emph{boundary at infinity} of $\Gamma$ (as defined by Gromov) is then realized as a $\Gamma$-invariant subset of the sphere at infinity of $X$.

In the past two decades, the development of the theory of \emph{Anosov groups} has built a nice framework to study geometric realizations of word-hyperbolic groups in higher rank Lie groups. Anosov representations were first introduced by Labourie in \cite{Labourie} for fundamental groups of closed negatively curved Riemannian manifolds, and the definition was later extended to more general word hyperbolic groups by Guichard--Wienhard in \cite{GW}. The definition was recently streamlined by various authors \cite{KLP1,BPS}, who proved in particular that a subgroup $\Gamma$ of a semisimple linear group $G$ satisfies a refinement of quasi-isometric embeddedness if and only if $\Gamma$ is Gromov hyperbolic and the inclusion is Anosov (see Definition \ref{d:Anosov}). 
 
Anosov subgroups of a semisimple linear Lie group $G$ have many good geometric and dynamical properties: they are quasi-isometrically embedded into the ambient group, are stable under small deformations, and their \emph{boundary at infinity} (as defined by Gromov) identifies with a compact invariant subset of some flag variety of $G$ (see \cite{Labourie, GW}). They are now commonly accepted as a good higher rank generalization of convex-cocompactness in rank 1. In fact, by the work of Danciger--Gu\'eritaud--Kassel \cite{DGK} and Zimmer \cite{Zimmer} every Anosov subgroup $\Gamma$ of $G$, after embedding $G$ in the appropriate linear group, acts convex-cocompactly on a stricly convex open domain of some projective space ${\bf P}(\mathbb{R}^d)$.
\medskip

Though not every Gromov hyperbolic group admits an Anosov representation, the only known obstruction so far seems to be non-linearity. Non-linear hyperbolic groups were first constructed by M. Kapovich \cite{Kap}: using Corlette's super-rigidity for lattices in $\Sp(k,1)$, he proved that the quotients of such lattices by a sufficiently large power of a non-trivial element are not linear. More recently, Canary, Stover and the second author constructed new examples by proving that sufficiently complicated amalgamated products of $\Sp(k,1)$-lattices are not linear. Our main result here is a new obstruction to realizing certain Gromov-hyperbolic groups as convex-cocompact groups. We prove the following:

\begin{theorem} \label{t:ExistenceLinearIndiscreteIntro}
There exists a Gromov-hyperbolic group which admits a faithful representation into $\Sp(k,1)$ but does not admit any discrete and faithful representation into any semisimple Lie group of real rank $1$.
\end{theorem}

Our example is the amalgamated product of two copies of a uniform $\Sp(k,1)$-lattice along a maximal infinite cyclic subgroup. In the process, we prove a general linearity result for the double of a torsion-free Anosov group along a maximal cyclic subgroup, as well as several non-linearity results for more complicated amalgamated products which simplify the construction of non-linear examples in \cite{CST}.

\subsection{Amalgamated products of Gromov hyperbolic groups}

Let $\Gamma_1$, $\Gamma_2$ and $W$ be groups and $i_1:W \rightarrow \Gamma_1$ and $i_2:W \rightarrow \Gamma_2$ be injective group homomorphisms. Recall that the \emph{amalgamated product} $\Gamma_1*_{i_1(W)=i_2(W)}\Gamma_2$ is the quotient of the free product $\Gamma_1*\Gamma_2$ by the normal subgroup generated by the set $\big \{g \big(i_1(w) i_2(w)^{-1} \big) g^{-1}: g \in \Gamma_1 \ast \Gamma_2, w\in W\big \}$. This operation is inspired by Van Kampen's theorem, which states that the fundamental group of a union of two open sets with connected intersection is the amalgamated product of the fundamental groups of the two open sets along the fundamental group of the intersection.

A particular case of amalgamated product is the \emph{double of $\Gamma$ along $W$}, denoted $\Gamma \ast_W\Gamma$, where $\Gamma_1 = \Gamma_2 = \Gamma$, $W$ is a subgroup of $\Gamma$ and $i_1$ and $i_2$ are both the inclusion of $W$ into $\Gamma$.

It is folklore knowledge that the free product of two word hyperbolic groups is again word hyperbolic. Bestvina--Feighn \cite{BF} have studied more generally the hyperbolicity of graphs of groups with word hyperbolic vertex and edge groups. In particular, they proved that the amalgamated product of hyperbolic groups along quasi-convex malnormal subgroups is hyperbolic (see also Theorem \ref{combination}). This applies for instance to amalgamations of torsion-free hyperbolic groups along maximal cyclic subgroups.

\subsection{Linearity and non-linearity}

Our first theorem is a linearity theorem for doubles of torsion free Anosov groups along a maximal cyclic subgroup:

\begin{theorem} \label{t:LinearityIntro}
Let $\Gamma$ be a torsion free Anosov subgroup of a semisimple linear group $G$ and $\langle w \rangle$ be a maximal cyclic subgroup of $\Gamma$. Then the double $\Gamma \ast_{\langle w \rangle}\Gamma$ of $\Gamma$ along $\langle w\rangle$ admits a faithful representation into $G$.
\end{theorem}

This result is somehow optimal. To stress this out, we will also prove several non-linearity results for more complicated amalgamated products.

\begin{theorem} \label{t:NonLinearAmalgamatedCyclicIntro}
Let $\Gamma_1$ and $\Gamma_2$ be lattices in $\Sp(k,1)$, $k \geq 2$, and $\langle w_i \rangle$ be a cyclic subgroup of $\Gamma_i$ for $i=1,2$. Assume that $w_1 \in \Gamma_1$ and $w_2 \in \Gamma_2$ have different translation lengths in the symmetric space of $\Sp(k,1)$. Then every linear representation of $\Gamma_1\ast_{w_1=w_2} \Gamma_2$ restricted on $\Gamma_1$ and $\Gamma_2$ has finite image.
\end{theorem}

\begin{corollary}
Let $\Gamma_1$ and $\Gamma_2$ be lattices in $\Sp(k,1)$, $k\geq 2$. Let $W$ be an infinite finitely generated group and $i_1$, $i_2$ be embeddings of $W$ into $\Gamma_1$ and $\Gamma_2$ respectively. Assume that there exists $w\in W$ such that $i_1(w)$ and $i_2(w)$ have different translation lengths. Then the amalgamated product $\Gamma_1 \ast_{i_1(W)=i_2(W)} \Gamma_2$ is not linear.
\end{corollary}

When $\Gamma_1$ and $\Gamma_2$ are torsion-free and $i_1(W)$ and $i_2(W)$ are maximal cyclic subgroups, this corollary gives new examples of non-linear hyperbolic groups by the Bestvina--Feighn combination theorem.

Theorem \ref{t:NonLinearAmalgamatedCyclicIntro} does not apply to doubles of a quaternionic lattice $\Gamma$ over a subgroup $W$. However, these also tend to be non-linear for a larger group $W$.

\begin{theorem} \label{t:NonLinearDoubleIntro}
Let $\Gamma$ be a uniform lattice in $\Sp(k,1)$, $k\geq 2$ and $W$ be a proper subgroup of $\Gamma$ which is not a uniform lattice in its Zariski closure. Then $\Gamma\ast_W\Gamma$ is not linear.
\end{theorem}

Again, this theorem can be applied to a quasi-convex malnormal free subgroups of $\Gamma$ (which exist by \cite[Thm. 6.7]{IKapovich}) and are not lattices in their Zariski closure in $\mathsf{Sp}(k,1)$, thus giving new constructions of non-linear hyperbolic groups. Both constructions are improvements on the main constructions of \cite{CST}.
\medskip

Finally, we point out the importance of the Anosov assumption in Theorem \ref{t:LinearityIntro} by proving the following:

\begin{theorem} \label{t:NonLinearitySLZ} For $n\geq 3$, there exist (many) maximal cyclic subgroups $\langle w \rangle$ of $\SL(n,\Z)$ for which the double of $\SL(n,\Z)$ along $\langle w \rangle$ is not linear.
\end{theorem}

\subsection{Discreteness and Anosov property}

The subgroup $\Gamma \ast_{\langle w \rangle}\Gamma$ of $G$ in Theorem \ref{t:LinearityIntro} has no reason to be discrete. In fact, it cannot be discrete if $\Gamma$ is already a uniform lattice in $G$. We will prove that some of them can never be embedded discretely into any rank $1$ Lie group.

\begin{restatable}{theorem}{indiscreteness} \label{t:IndiscretenessIntro}
Let $\Gamma$ be a uniform lattice in $\Sp(k,1)$, $k \geq 4$, and $\langle w \rangle$ be an infinite maximal cyclic subgroup of $\Gamma$. Then the group $\Gamma\ast_{\langle w \rangle} \Gamma$ does not admit a discrete and faithful representation into any semisimple Lie group of rank $1$. \end{restatable}

Recall that if $\Gamma$ is torsion-free then $\Gamma\ast_{\langle w \rangle} \Gamma$ is word hyperbolic by the Bestvina--Feighn combination theorem (Theorem \ref{combination}) and isomorphic to a (dense) subgroup of $\Sp(k,1)$ by Theorem \ref{t:LinearityIntro}. The group $\Gamma \ast_{\langle w \rangle} \Gamma$ is an example satisfying the conclusion of Theorem \ref{t:ExistenceLinearIndiscreteIntro}. To our knowledge, this is the first example of a linear word-hyperbolic group which is not virtually isomorphic to a convex cocompact group of a \hbox{simple Lie group of real rank $1$.}

Note that Theorem \ref{t:IndiscretenessIntro} contrasts with the following theorem of Baker--Cooper in real hyperbolic geometry: 

\begin{theorem}[Baker--Cooper \cite{Baker-Cooper}]
Let $\Gamma$ be a convex-cocompact group of isometries of the real hyperbolic space $\mathbb{R}{\bf H}^k$ and $\langle w \rangle$ be an infinite maximal cyclic subgroup of $\Gamma$. Then there exists a finite index subgroup $\Gamma'$ of $\Gamma$ containing $\langle w \rangle$ such that $\Gamma'\ast_{\langle w \rangle} \Gamma'$ admits a convex-cocompact representation into $\Isom(\mathbb{R}{\bf H}^{2k-1})$.
\end{theorem}

\subsection{Further questions and perspectives}

The present work leaves open the following question:

\begin{ques}\label{ques:AmalgamatesAnosov}
Let $\Gamma$ be an Anosov subgroup of $G$ and $\langle w \rangle$ a maximal cyclic subgroup of $\Gamma$. Does $\Gamma\ast_{\langle w\rangle } \Gamma$ admit an Anosov representation \textup{(}possibly in some larger group\textup{)}?
\end{ques}

We strongly believe in the following weaker statement, motivated partly by Baker--Cooper's theorem above:

\begin{conj}\label{conj:AmalgamatesAnosov}
Let $\Gamma$ be an Anosov group and $\langle w \rangle$ a maximal cyclic subgroup of $\Gamma$. Then there exists a finite index subgroup $\Gamma'$ of $\Gamma$ containing $w$ such that $\Gamma'\ast_{\langle w \rangle} \Gamma'$ admits an Anosov representation.
\end{conj}

\begin{remark} Perhaps it is worth noting that, even though $\Gamma'$ is a proper finite index subgroup of $\Gamma$ with $w \in \Gamma'$, the group $\Gamma'\ast_{\langle w \rangle} \Gamma'$ has infinite index in $\Gamma \ast_{\langle w \rangle} \Gamma$. Hence Conjecture \ref{conj:AmalgamatesAnosov} does not readily answer to Question \ref{ques:AmalgamatesAnosov}. \end{remark}

Note that the analoguous question for free products has already been studied. Dey--Kapovich--Leeb in \cite{Dey-Kapovich-Leeb} and Dey--Kapovich \cite{DK} have studied when the group generated by two $P$-Anosov subgroups of a group $G$ is a $P$-Anosov free product into $G$. Danciger--Gu\'eritaud--Kassel have also announced (see \cite[Prop. 12.5]{DGK}) that the free product of two infinite Anosov subgroups of $\mathsf{PGL}(d,\mathbb{R})$ also admits an Anosov representation into $\mathsf{PGL}(m,\mathbb{R})$ for some $m \geq d$. The particular case of convex-cocompact groups in rank 1 Lie groups seems to be folklore.

Let us finally mention that the work of Agol and Wise give many examples of hyperbolic groups admitting discrete and faithful linear representations: Wise in \cite{Wise} proved that $C'(\frac{1}{6})$-small cancellation groups are cubulated, and Agol \cite{Agol} proved (relying also on the work of Haglund--Wise \cite{Haglund-Wise}) that cubulated hyperbolic groups are virtually special and embed into some $\GL(d,\Z)$.
\medskip

These various results raise the following general question:

\begin{ques} Is there a linear hyperbolic group that does not admit a discrete and faithful representation in any (real) linear group ?
\end{ques}

\subsection{Strategy of the proofs}

We will derive Theorem  \ref{t:LinearityIntro} from the linearity for the HNN extension $\Gamma\ast_{\langle w \rangle} \equaldef \left \langle \Gamma, t \mid twt^{-1}=w \right \rangle$:

\begin{theorem} \label{t:LinearityBis}
Let $\Gamma$ be a torsion free Anosov subgroup of a semisimple linear group $G$ and $\langle w \rangle$ a maximal cyclic subgroup of $\Gamma$. Then there exists $t\in G$ which commutes with $w \in \Gamma$ such that the subgroup of $G$ generated by $\Gamma$ and $t$ is isomorphic to the HNN extension $\Gamma\ast_{\langle w \rangle}$. 
\end{theorem}

The element $t \in G$ then satisfies the requirements of Theorem \ref{t:LinearityIntro}. In fact, the group $\Gamma\ast_{\langle w \rangle}$ has a morphism onto $\Z$, the kernel of which is spanned by all the $t^i \Gamma t^{-i}, i\in \Z$. This kernel is isomorphic  to the amalgamated product of infinitely many copies of $\Gamma$ over $\langle w \rangle$.

The other negative theorems (Theorems \ref{t:NonLinearAmalgamatedCyclicIntro}, \ref{t:NonLinearDoubleIntro}, \ref{t:NonLinearitySLZ} and \ref{t:IndiscretenessIntro}) are all based on the same strategy as in \cite{CST}: A linear representation of an amalgamated product $\Gamma_1\ast_W \Gamma_2$ is given by two representations $\rho_1$ and $\rho_2$ of $\Gamma_1$ and $\Gamma_2$ respectively which coincide on $W$. We consider situations where the superrigidity theorems of Margulis \cite{Mar} and Corlette \cite{Corlette} give such strong constraints on $\rho_1$ and $\rho_2$ so that adding the compatibility condition on $W$ leads to a contradiction.

The paper is organized as follows. In Section \ref{s:background}, we recall some elements on the structure of semisimple linear groups and the main properties of their Anosov subgroups. In Section \ref{s:linearity} we prove our linearity theorem, Theorem ~\ref{t:LinearityBis}, from which Theorem \ref{t:LinearityIntro} follows. In Section \ref{s:superrigidity} we recall the precise statements of Margulis' and Corlette's superrigidity theorems, which we use in Section \ref{s:proofs} to prove the remaining results. 

\subsection*{Acknowledgements} We would like to thank Richard Canary, Beatrice Pozzetti and Matthew Stover for helpful conversations, as well as Olivier Benoist for his Galois theoretic help with the proof of Lemma \ref{l:ZariskiColsureCyclicSL(d,Z)}. This project received funding from the European Research Council (ERC) under the European's Union Horizon 2020 research and innovation programme (ERC starting grant DiGGeS, grant agreement No 715982).
\medskip

\section{Background} \label{s:background}

In this section, we provide some background on amalgamated products and HNN extensions, Lie theory and Anosov representations.

\subsection{Amalgamated products and HNN extensions} \label{ss:Amalgamated products and HNN}
We refer the reader to \cite[Ch. IV]{LS} and \cite{Serre} for more background on amalgamated products and HNN extensions. For a group $H$ and a subset $S$ of $H$, we denote by $\langle \langle S \rangle \rangle$ the normal subgroup of $H$ generated by $\big\{hsh^{-1}:s \in S, h \in H\big\}$. 

Let $\Gamma_1$, $\Gamma_2$ and $W$ be three groups and $\iota_1: W \to \Gamma_1$ and $\iota_2:W\to \Gamma_2$ two injective\footnote{The general definition does not require injectivity, but this assumption will simplify the exposition here, and is enough for our purposes.} morphisms.

\begin{definition}
The amalgamated product of $\Gamma_1$ and $\Gamma_2$ along $W$ is the group
\[\Gamma_1 \ast_{\iota_1 = \iota_2} \Gamma_2 = \Gamma_1\ast \Gamma_2 / \langle \langle \iota_1(w)\iota_2^{-1}(w), w\in W \rangle \rangle~.\]
\end{definition}

Let $T_i$ ($i=1,2$) be a sets of right coset representatives of $\iota_i(W)$ in $\Gamma_i$. A normal form is a sequence $(w_0,\ldots, w_n)$, $n \geq 0$, with the following properties: $w_0 \in W$, if $n \geq 1$, $w_i \in T_1\smallsetminus \{1\}$ or $w_i \in T_2 \smallsetminus \{1\}$ for $1 \leq i \leq n$ and for every $1\leq i\leq n-1$ we have $x_i \in T_1$ and $x_{i+1}\in T_2$ or $x_i \in T_2$ and $x_{i+1}\in T_1$. Every element $g \in \Gamma \ast_{\iota_1 = \iota_2}B_1$ has a unique representation $g=\iota_1(w_0) w_1 \cdots w_n$ for some normal form $(w_0,\ldots , w_n)$ (see \cite[Ch. IV, Thm. 2.6]{LS}). In particular, $\Gamma_1$ and $\Gamma_2$ naturally identify to subgroups of $\Gamma$. 

The amalgamated product $\Gamma_1 \ast_{\iota_1=\iota_2} \Gamma_2$ satisfies the following universal property: 

\begin{proposition}
For any group $G$ and any homomorphisms $\rho_1:\Gamma_1\to G$ and $\rho_2:\Gamma_2 \to G$ such that $\rho_1\circ i_1 = \rho_2\circ i_2$, there exists a unique homomorhism $\rho:\Gamma_1 \ast_{\iota_1=\iota_2} \Gamma_2 \to G$ whose restriction to $\Gamma_i$ is $\rho_i$ for $i=1,2$. 
\end{proposition}

If the morphisms $\iota_i$ are implicit, we will sometimes denote the amalgamated product by $\Gamma_1 \ast_{W_1=W_2} \Gamma_2$, where $W_i = \iota_i(W)$, or even $\Gamma_1\ast_W \Gamma_2$ when this does not bring any confusion. When $\Gamma_1 = \Gamma_2= \Gamma$, $W$ is a subgroup of $\Gamma$ and $\iota_1, \iota_2$ are the inclusion, we call $\Gamma \ast_W \Gamma$ the \emph{double} of $\Gamma$ along $W$.

The following fact, which will be useful later, is a straightforward consequence of the uniqueness of normal forms in amalgamated free products.

\begin{fact} \label{nonfaithful} Let $\Gamma$ be a group, $\Gamma_1$ and $\Gamma_2$ be subgroups of $\Gamma$ such that $\Gamma_1$ is a proper subgroup of $\Gamma_2$. Then the natural group homomorphism $\pi:\Gamma \ast_{\Gamma_1} \Gamma \rightarrow \Gamma \ast_{\Gamma_2} \Gamma$ is not injective. \end{fact}
\medskip

Let now $\Gamma$ be a group, $W$ a subgroup of $\Gamma$ and $\varphi:W \to \Gamma$ a morphism.

\begin{definition}
The HNN extension of $\Gamma$ {\em relative to $\varphi$} is the group $$\Gamma\ast_{\varphi}=\Gamma \ast \langle t \rangle \big/ \langle \langle twt^{-1}\varphi(w)^{-1},a \in W\rangle \rangle.$$
\end{definition}

Britton's lemma (see \cite[Ch. IV, \S2]{LS}) states that a word $g\in \Gamma \ast_{\varphi}$ is non-trivial if $g=g_0t^{m_0}g_1t^{m_1}\cdots g_s t^{m_s}$, where $g_0,\ldots, g_s \in \Gamma$, $m_1,\ldots, m_s\in \{-1,1\}$ and there are no subwords of the form $tg_it^{-1}$ where $g \in W\smallsetminus \{1\}$ and $t^{-1} g_i t$ where $g \in \varphi(W)\smallsetminus\{1\}$. In particular, $\Gamma$ embeds as a subgroup of $\Gamma\ast_\varphi$ (see \cite[Ch. IV, Thm. 2.1]{LS}).

The HNN extension $\Gamma \ast_\phi$ satisfies a universal property:
\begin{proposition}\label{univpropHNN}
Let $G$ be a group, $\rho:\Gamma \to G$ a homomorphism and $h\in G$ such that $\rho(\varphi(w)) = h\rho(w)h^{-1}$ for all $w\in W$. Then there exists a unique homomorphism $\rho_h: \Gamma \ast_\varphi \to G$ such that ${\rho_h}_{\vert \Gamma} = \rho$ and $\rho_h(t)=h$.
\end{proposition}

When $\varphi$ is the inclusion of $W$ in $\Gamma$, we simply denote the HNN extension by $\Gamma \ast_W$. This HNN extension admits a surjective group homomorphism $\pi:\Gamma\ast_{W}\rightarrow \mathbb{Z}$ mapping $\Gamma$ to $0$ and $t$ to $1$. Its kernel is the normal subgroup $\langle \langle \Gamma \rangle \rangle=\big \langle t^{i}at^{-i}: g \in \Gamma, i \in \mathbb{Z}\big \rangle$, which is isomorphic to the amalgamated product of countably many copies of $\Gamma$ along $W$.\\

A subgroup $W$ of $\Gamma$ is \emph{malnormal} if, for every $g\in \Gamma \smallsetminus W$, we have $gWg^{-1} \cap W =\{1\}$ . When $\Gamma$ is word-hyperbolic, $W$ is \emph{quasi-convex} in $\Gamma$ if and only if the inclusion $W \hookrightarrow \Gamma$ induces  quasi-isometric embedding of their Cayley graphs. The \emph{Bestvina--Feighn combination theorem} asserts (in particular), that amalgamations of hyperbolic groups along malnormal quasi-convex subgroups are hyperbolic.

\begin{theorem}[Bestvina--Feighn \cite{BF}] \label{combination}
Let $\Gamma_1$, $\Gamma_2$ and $W$ be word-hyperbolic groups and $\iota_1 : W \hookrightarrow \Gamma_1$ and $\iota_2:W \hookrightarrow \Gamma_2$ be quasi-convex embeddings such that $\iota_1(W)$ is malnormal in $\Gamma_1$. Then $\Gamma_1\ast_{\iota_1=\iota_2} \Gamma_2$ is word hyperbolic.
\end{theorem}

Note that a maximal cyclic subgroup\footnote{i.e. is not contained in a larger cyclic subgroup.} of a torsion-free word-hyperbolic group is  malnormal. Kapovich \cite{IKapovich} proved that every non-elementary word-hyperbolic group contains malnormal quasi-convex free subgroups with $2$ generators.

\subsection{Lie theory} Let us fix some notation, mainly following \cite[\S 3.2]{GW}. Throughout this paper, we consider $G$ a linear, non-compact, semisimple real algebraic Lie group and denote by $\mathfrak{g}$ its Lie algebra. Let $\textup{Ad}:G \rightarrow \mathsf{GL}(\mathfrak{g})$ and $\textup{ad}: \mathfrak g \to \mathrm{End}(\mathfrak g)$ denote the adjoint representations of $G$ and $\mathfrak g$, and let $\exp:\mathfrak{g}\rightarrow G$ be the exponential map. The Killing form $B: \mathfrak{g}\times \mathfrak{g} \rightarrow \mathbb{R}$ is the bilinear form $B(X,Y)=\textup{tr}(\textup{ad}_{X}\textup{ad}_{Y})$. It is invariant under the adjoint action and is non-degenerate as soon as $\mathfrak{g}$ is semisimple. Let us fix $K$ a maximal compact subgroup of $G$, unique up to conjugation. We have the associated decomposition $$\mathfrak{g}=\mathfrak{k}\oplus \mathfrak{p}$$ where $\mathfrak{k}=\textup{Lie}(K)$ and $\mathfrak{p}$ is its orthogonal complement with respect to $B$. We choose a Cartan subspace $\mathfrak{a}$, i.e. a maximal abelian subalgebra of $\mathfrak{g}$ contained in $\mathfrak{p}$. The {\em real rank} of $G$ is the dimension of $\mathfrak{a}$ as a real vector space. 
 
\par The co-diagonalization of the adjoint action of $\mathfrak{a}$ decomposes $\mathfrak g$ as the direct sum of \emph{root spaces}:$$\mathfrak{g}=\mathfrak{g}_0\oplus \bigoplus_{\alpha \in \Sigma} \mathfrak{g}_{\alpha}~,$$ where $\mathfrak{g}_{\alpha}=\{X \in \mathfrak{g}:\mathsf{ad}_H(X)=\alpha(H)X \ \forall H \in \mathfrak{a}\}$ and $\Sigma=\{\alpha \in \mathfrak{a}^{\ast}:\mathfrak{g}_{\alpha}\neq 0\}$ is the set of \emph{restricted roots}. After chosing a vector $u \in \mathfrak a \smallsetminus \bigcup_{\alpha \in \Sigma} \ker \alpha$, we define the set of \emph{positive roots}
\[\Sigma^+ = \{\alpha \in \Sigma \mid \alpha(u) >0\}\]
and the \emph{dominant Weyl chamber}
\[\mathfrak a^+ =\left \{H \in \mathfrak{a}:\alpha(H)\geq 0 \ \forall \alpha \in \Sigma^{+} \right \}~.\]
Finally, a positive root is \emph{simple} if it cannot be written as the sum of two positive roots. The set of simple roots $\Delta$ is a basis of $\mathfrak a^*$ and the dominant Weyl chamber is the associated positive quadrant.

The {\em Cartan decomposition} writes every element $g\in G$ in the form $g=k\exp(\mu(g))k'$ for $k,k' \in K$ and $\mu(g)\in \mathfrak a^+$. The vector $\mu(g)$ is unique and called the {\em Cartan projection} of $G$. The {\em Jordan projection} $\overrightarrow{\ell}:G \rightarrow \mathfrak{a}^{+}$ can be defined by \[\overrightarrow{\ell}(g)=\lim_{n \rightarrow \infty} \frac{\mu(g^n)}{n}~.\]

\par Let $\theta \subset \Delta$ be a subset of simple restricted roots. Set $$\mathfrak{a}_{\theta}\equaldef \bigcap_{\alpha \in \Delta \smallsetminus \theta}\textup{ker}(\alpha)$$ and let $Z_{K}(\mathfrak{a}_{\theta})$ be the centralizer of $\mathfrak{a}_{\theta}$ in $K$. One associates to $\theta$ a pair of \emph{opposite parabolic subgroups} $(P_{\theta}^{+},P_{\theta}^{-})$ defined by:
\[ P_{\theta}^{\pm} =Z_K (\mathfrak{a}_{\theta})\exp(\mathfrak{a}) \exp \Big(\bigoplus_{\alpha \in \Sigma^{+}} g_{\pm \alpha}\Big ).\]
The subgroup $L_{\theta} \equaldef P_{\theta}^{+}\cap P_{\theta}^{-}$ is a Levi factor of both $P_\theta^+$ and $P_\theta^-$, which admits the Cartan decomposition $$L_{\theta}=Z_{K}(\mathfrak{a}_{\theta})\exp\big( \mathfrak{a}_{L_{\theta}}^{+} \big)Z_{K}(\mathfrak{a}_{\theta})$$ where $\mathfrak{a}_{L_{\theta}}^{+}=\big \{ H \in \mathfrak{a}: \alpha(H)\geq 0, \ \forall \alpha \in \theta\big\}$.

An element $g\in G$ is called {\em $\theta$-proximal} if $\alpha(\lambda(g))>0$ for every $\alpha \in \theta$. For a group $\mathsf{H}$, a representation $\rho:\mathsf{H}\rightarrow G$ is calles $\theta$-proximal if $\rho(\mathsf{H})$ contains a $\theta$-proximal element.

\subsection{The example of $\SL(d,\Real)$} 
We endow $\mathbb{R}^d$ with its standard inner product and denote by $(e_1,\ldots,e_d)$ its canonical othonormal basis. For $J \subset \{1,\ldots,d\}$ we set $\langle \{e_j: j \in J\}\rangle^{\perp}\equaldef \langle \{e_i:i \notin J\} \rangle$.

A standard choice of compact subgroup of $\SL(d,\Real)$ is the special orthogonal group $\mathsf{SO}(d)=\big \{g \in \SL(d,\mathbb{R}):gg^{t}=\mathrm I_d \big \}$. One can chose as a Cartan subspace the space $\mathfrak{a}=\mathsf{diag}_{0}(d)$ of diagonal matrices with trace zero, and as a dominant Weyl chamber the cone of traceless diagonal matrices with diagonal coefficients in non-increasing order. The restricted roots are the forms $\varepsilon_i - \varepsilon_j$ where $\varepsilon_i \in \mathfrak{a}^{\ast}$ is the projection to the $(i,i)$ entry, and the root space associated to $\varepsilon_i - \varepsilon_j$ (for $i\neq j$) is $\mathbb{R}E_{ij}$, where $E_{ij}$ is the $d \times d$ elementary matrix with $1$ at the $(i,j)$ entry and $0$ everywhere else. The root $\varepsilon_i - \varepsilon_j$ is positive when $i>j$ and simple when $j=i+1$. 

Given $g \in \mathsf{SL}(d, \mathbb{R})$, denote by $\lambda_1(g) \geq \lambda_2(g)\geq \ldots \geq \lambda_d(g)$ the moduli of the eigenvalues of $g$ in decreasing order (counting multiplicity), and by $\sigma_i(g)$ the $i$\textsuperscript{th} singular value of $g$, defined by the relation $\sigma_i(g)=\sqrt{\lambda_i(g^tg)}$. The Cartan and Jordan projections of $g \in \mathsf{SL}(d,\mathbb{R})$ \hbox{are given by the vectors} \begin{align*} \mu(g) &=\big(\log \sigma_1(g),\ldots,\log \sigma_d(g)\big)\\ \overrightarrow{\ell}(g) & =\big(\log \lambda_1(g),\ldots,\log\lambda_d(g)\big)\end{align*} respectively. For $1\leq i \leq d-1$, the matrix $g$ is called {\em $i$-proximal} if $\lambda_i(g)>\lambda_{i+1}(g)$.

Now, for $\theta =\{ \varepsilon_{i_1}- \varepsilon_{i_1+1}, \ldots ,\varepsilon_{i_k}- \varepsilon_{i_k+1}\} \subset \Delta$, the associated parabolic subgroups $P_\theta^+$ and $P_\theta^-$ are respectively the stabilizers of the flag
$$\langle e_1,\ldots, e_{i_1}\rangle \subset \langle e_1,\ldots , e_{i_2} \rangle \subset \ldots \subset \langle e_1,\ldots, e_{i_k}\rangle$$ and the stabilizer of the flag $$\langle e_1,\ldots, e_{i_k}\rangle^{\perp} \subset \langle e_1,\ldots , e_{i_2}^\perp \rangle \subset  \ldots \subset \langle e_1,\ldots, e_{i_1}\rangle^{\perp}.$$

\subsection{Anosov representations.} 

More recently, Kapovich--Leeb--Porti in \cite{KLP1} and Bochi--Potrie--Sambarino in \cite{BPS}, characterized Anosov representations into $G$ entirely in terms of their Cartan projections. Moreover, this particularly synthetic characterization does not a priori assume hyperbolicity of the domain group. We use it here as a definition.

Let $\Gamma$ be a finitely generated group. We fix a left invariant word metric on $\Gamma$ and denote by $|\gamma|$ the distance of and element $\gamma\in \Gamma$ from the identity element $e \in \Gamma$.

\begin{definition}[\cite{KLP1, BPS}] \label{d:Anosov} Let $\theta \subset \Delta$ be a \textup{(}non-empty\textup{)} subset of simple restricted roots of $G$. A representation $\rho:\Gamma \rightarrow G$ is called \emph{$P_{\theta}$-Anosov} or \emph{$\theta$-Anosov} if there exist constants $C,c>0$ such that $$\alpha \left(\mu(\rho(\gamma)) \right) \geq c|\gamma|-C$$ for every $\alpha \in \theta$ and all $\gamma \in \Gamma$.
\end{definition}

If $\Gamma$ is a finitely generated subgroup of $G$, we call $\Gamma$ \emph{$P_{\theta}$-Anosov} or \emph{$\theta$-Anosov} if the inclusion $\Gamma \hookrightarrow G$ is $\theta$-Anosov. For an abstract group $\Gamma$, any $\theta$-Anosov representation $\rho:\Gamma \to G$ has finite kernel and its image is a $\theta$-Anosov subgroup of $G$.

Anosov groups enjoy many nice geometric and dynamical properties. In particular, they are word-hyperbolic and their Gromov boundary identifies with a compact subset of some flag variety of $G$. To be more precise, let $\Gamma$ be a word-hyperbolic group and $\partial_{\infty}\Gamma$ its Gromov boundary. Every infinite order element $\gamma \in \Gamma$ has exactly two fixed points on $\partial_{\infty}\Gamma$ denoted by $\gamma^{+}$ and $\gamma^{-}$ called respectively the attracting and repelling fixed point of $\gamma$.

\begin{theorem}[Labourie \cite{Labourie}, Guichard--Wienhard \cite{GW}, Kapovich--Leeb--Porti \cite{KLP1}] \label{t:BoundaryMaps}
Let $\rho:\Gamma \to G$ be a $\theta$-Anosov representation. Then there exists a unique pair of continuous, $\rho$-equivariant injective maps $$\big(\xi_{\rho}, \xi_{\rho}^{-}\big):\partial_{\infty}\Gamma \rightarrow G/P_{\theta}^{+} \times G/P_{\theta}^{-}$$ called the {\em Anosov limit maps} of $\rho$, satisfying the following properties:
\medskip
\begin{itemize}
\item $\xi_{\rho}$ and $\xi_{\rho}^{-}$ are transverse: for every pair $(x,y)$ of distinct points of $\partial_{\infty}\Gamma$, there exists $h\in G$ such that $\xi_{\rho}(x)=hP_{\theta}^{+}$ and $\xi_{\rho}^{-}(y)=hP_{\theta}^{-}$.

\item The maps $\xi_{\rho}$ and $\xi_{\rho}^{-}$ are dynamics preserving: for every infinite order element  $\gamma \in \Gamma$, $\rho(\gamma)$ is $\theta$-proximal and the points $\xi_\rho(\gamma^{+})$ and $\xi_\rho^-(\gamma^{+})$ are the attracting fixed points of $\rho(\gamma)$ in $G/P_\theta^+$ and $G/P_\theta^-$ respectively. 
\end{itemize}

\end{theorem}

For more background on Anosov representations we refer to the survey paper \cite{Canary-notes}.

\subsection{The rank 1 case.}

When the linear group $G$ has rank $1$, its symmetric space $G/K$ is Gromov hyperbolic and its boundary at infinity coincides with $G/P$, where $P$ is the unique proper parabolic subgroup of $G$ up to conjugation.

A representation $\rho: \Gamma \to G$ is $P$-Anosov if and only if it is a quasi-isometric embedding, and the existence of an Anosov limit map $\xi_\rho : \partial_\infty \Gamma \to G/P$ reduces to a property of quasi-isometric embeddings between Gromov hyperbolic spaces.

Finally, $\rho$ acts properly discontinuously and cocompactly on the convex hull of $\xi_\rho(\partial_\infty \Gamma)$ in $G/K$. Thus, in rank $1$, the notion of Anosov representation coincides with the classical notion of \emph{convex cocompact representation}.

\subsection{Projective Anosov representations}

A particular case of Anosov property is the \emph{projective Anosov} property of representations into $\SL(d,\Real)$.

\begin{definition}
A linear representation $\rho:\Gamma \rightarrow \mathsf{SL}(d,\mathbb{R})$ is called {\em projective Anosov} if there exist constants $C,a>0$ with $$\frac{\sigma_1(\rho(\gamma))}{\sigma_2(\rho(\gamma))}\geq Ce^{a|\gamma|_{\Gamma}}$$ for every $\gamma \in \Gamma$. 
\end{definition}

In other words, a projective Anosov representation is a $\theta$-Anosov representation for
\[\theta= \{\varepsilon_1-\varepsilon_2\}~.\]
Such a representation induces boundary maps $\xi_\rho$ and $\xi_\rho^-$ with values into $\mathbf{P}(\mathbb{R}^d)$ and $\mathsf{Gr}_{d-1}(\mathbb{R}^d)$ respectively, and the transversality condition means for all $x\neq y \in \partial_\infty \Gamma$, 
\[\xi_\rho(x)\notin \ker \xi_\rho^-(y)~.\]

\subsection{Restricted weights of representations of $G$} Let $G$ be a linear real semisimple Lie group, $\tau:G \rightarrow \mathsf{GL}(d,\mathbb{R})$ be an irreducible linear representation of $G$ and fix $\mathfrak{a} \subset \mathfrak{g}$ a Cartan subspace of $G$. Up to conjugation, we may assume that $d\tau (\mathfrak{a})$ is contained in $\mathsf{diag}_0(d)$ and hence there exists a decomposition of $\mathbb{R}^d$ into weight spaces: $$\mathbb{R}^d=V^{\chi_1}\oplus \cdot \cdot \cdot \oplus V^{\chi_{d}}$$ where $\chi_1,\ldots,\chi_d \in \mathfrak{a}^{\ast}$ and $V^{\chi_{i}} \equaldef \big \{v \in \mathbb{R}^d:d\tau(H)v=\chi_{i}(H)v, \ \forall H \in \mathfrak{a}\big \}$ for each $1 \leq i \leq d$. The linear forms $\chi_1,\ldots, \chi_{d} \in \mathfrak{a}^{\ast}$ are called the {\em restricted weights of} $\tau$ and are integral combinations of the fundamental weights $\{\omega_{\alpha}\}_{\alpha \in \Delta}$ of $G$ with respect to the set of simple roots $\Delta$. There is also a partial order on the set of distinct restricted weights of $\tau$: given two weights $\chi \neq \chi'$, $\chi>\chi'$ if and only if $\chi-\chi' \in \sum_{\alpha \in \Delta}\mathbb{R}^{+} \alpha$. For this partial order, among the distinct weights of $\tau$ there exists a unique maximal element, denoted by $\chi_{\tau}$, \hbox{called the {\em highest weight} of $\tau$.}

In fact, one can turn any Anosov representation into some Lie group $G$ into a projective Anosov representation after composing with a suitable linear representation of $G$ (see for instance \cite[Prop. 4.3]{GW}, \cite[Prop. 3.5]{GGKW}):

\begin{proposition} \label{higherdimension} Let $G$ a real semisimple Lie group and $\theta \subset \Delta$ a \textup{(}non-empty\textup{)} subset of simple restricted roots of $G$. There exists $d=d(G,\theta)$ and an irreducible $\theta$-proximal representation $\tau:G \rightarrow \mathsf{GL}(d, \mathbb{R})$ with the following properties:

\begin{itemize}
\item[(i)] A representation $\rho:\Gamma \rightarrow G$ of a word hyperbolic group $\Gamma$ is $\theta$-Anosov if and only if $\tau\circ \rho$ is projective Anosov.

\item[(ii)] Let $\chi \neq \chi_{\tau}$ be a restricted weight of $\tau$. If $H \in \overline{\mathfrak{a}}^{+}$ with $\alpha(H)>0$ for every $\alpha \in \theta$, then $\chi_{\tau}(H)>\chi(H)$. Moreover, if $\alpha(H)=0$ for every $\alpha  \in \Delta \smallsetminus \theta$ and $\alpha(H)\in \mathbb{N}^{\ast}$ for every $\alpha \in \theta$, then $\chi_{\tau}(H)-\chi(H)\in \mathbb{N}^{\ast}$. 
\end{itemize}

\end{proposition}

\begin{remark}\label{Plucker} \normalfont{(see \cite[Prop. 3.3 \& Prop. 3.5(i)]{GGKW}) For a $\theta$-Anosov representation \hbox{$\rho:\Gamma \rightarrow G$}, the Anosov limit maps of the projective Anosov representation $\tau\circ \rho:\Gamma \rightarrow \mathsf{GL}(d,\mathbb{R})$ can be obtained from $\tau$ and the limit maps $(\xi_{\rho},\xi_{\rho}^{-})$ of $\rho$ as follows. Let $V^{\chi_{\tau}}$ be the one dimesnional weight space corresponding to the highest weight and $V^{<\chi_{\tau}}$ be the direct sum of the weight spaces different from $V^{\chi_{\tau}}$. There exist $\tau$-equivariant embeddings $\iota^{+}:G/P_{\theta}^{+} \xhookrightarrow{} \mathbb{P}(\mathbb{R}^d)$ and $\iota^{-}:G/P_{\theta}^{-}\xhookrightarrow{} \mathsf{Gr}_{d-1}(\mathbb{R}^d)$ defined as follows: $$\iota_{\tau}^{+}\big(gP_{\theta}^{+}\big)=\big[\tau(g)V^{\chi_{\tau}} \big] \ \ \textup{and} \ \ \iota_{\tau}^{-}\big(gP_{\theta}^{-}\big)=\tau(g)V^{<\chi_{\tau}}$$ for $g \in G$. Then the pair of the Anosov limit maps of $\tau \circ \rho$ is $\big(\iota_{\tau}^{+}\circ \xi_{\rho}, \iota_{\tau}^{-}\circ \xi_{\rho}^{-}\big)$.
}\end{remark}

\subsection{Quaternionic hyperbolic spaces and groups.} 

Let $\mathbb{H}=\mathbb{R}\oplus \mathbb{R}i\oplus \mathbb{R}j\oplus \mathbb{R}k$ be Hamilton's quaternion algebra. For $z=a+bi+cj+dk \in \mathbb{H}$, denote by $\overline{z}=a-bi-cj-dk$ the conjugate of $z$. For a matrix $g=(g_{ij})_{ij} \in \textup{Mat}_{n \times n}(\mathbb{H})$, $g^{\ast}=(\overline{g_{ji}})_{ij}$ is the conjugate transpose of $g$. For $m \geq 1$, let $J_{m}=\textup{diag}\big(1,\ldots,1,-1\big)$. The projectivization ${\bf P}(\mathbb{H}^{m+1})$ is the set of equivalence classes of vectors in $\mathbb{H}^{m+1}$, where $u,v \in \mathbb{H}^{m+1}$ are equivalent if $u=vz$ for some $z \in \mathbb{H}-\{0\}$. The quaternionic hyperbolic space $\mathbb{H}{\bf H}^m$ is the open subset of the projective space $\mathbf{P}(\mathbb{H}^m)$ given by $$\mathbb{H}{\bf H}^m=\big\{\big[z_0:\ldots:z_{m-1}:z_m\big] \in {\bf P}(\mathbb{H}^{m+1}):|z_0|^2+\ldots+|z_{m-1}|^2<|z_m|^2\big \}$$ We consider the symplectic unitary group $\mathsf{Sp}(m,1)=\big\{ g \in \mathsf{GL}(m+1,\mathbb{H}): g^{\ast}J_mg=J_m \big\}$ and the compact subgroup $\mathsf{Sp}(m)=\big \{g \in \mathsf{GL}(m,\mathbb{H}):gg^{\ast}=I_m\big \}$. The Lie group $\mathsf{Sp}(m,1)$ preserves and acts transitively on $\mathbb{H}{\bf H}^m$. The stabilizer of ${\bf v_0}=\big[\{0\}^{m}:1\big]$ in $\mathsf{Sp}(m,1)$ is the group $\mathsf{Sp}(m)\times \mathsf{Sp}(1)$ which is also the, unique up to conjugation, maximal compact subgroup  of $\mathsf{Sp}(m,1)$.
\medskip

We shall use the following fact, which is almost immediate from the classification of the totally geodesic subspaces of $\mathbb{H}{\bf H}^m$.
\medskip

\begin{fact}\label{fact} Let $m \geq 1$ and $L$ be a geodesic in $\mathbb{H}{\bf H}^m$. There exists a unique quaternionic line $L'$ in $\mathbb{H}{\bf H}^m$ containing $L$. \end{fact}

\begin{proof} Up to translating $L$ by an element $g \in \mathsf{Sp}(m,1)$, we may assume that $$L=\big \{\big[\{0\}^{m-1}:\tanh t:1\big]: t \in \mathbb{R}\big \}=\mathbb{P}\big(\{0\}^{m-1}\times \mathbb{R}^2\big)\cap \mathbb{H}{\bf H}^m$$ Obviously, the quaternionic line $L'={\bf P}\big(\{0\}^{m-1}\times \mathbb{H}^2\big)\cap \mathbb{H}{\bf H}^m$ contains $L$. Suppose that $P$ is a quaternionic geodesic subspace of $\mathbb{H}{\bf H}^m$ containing $L$. By the classification of totally geodesic subspaces of $\mathbb{H}{\bf H}^m$ (see for example \cite[Thm. 2.12]{M}), there exists a real subspace $W$ of $\mathbb{H}^{m+1}$, right invariant by $\mathbb{H}$, so that $P={\bf P}(W)\cap \mathbb{H}{\bf H}^m$. Since $L$ is contained in $P$, $W$ contains $\{0\}\times \mathbb{R}^2$ and so contains $\{0\}\times \mathbb{H}^2$. It follows that $P$ contains $L'$. Therefore, $L'$ is unique.\end{proof}

\section{Linear groups from Anosov representations} \label{s:linearity}

Following an idea of Shalen (see \cite[Thm. 2]{Shalen}), we prove that certain amalgamated free products and HNN extensions of groups admitting Anosov representations over cyclic subgroups are linear. Using similar methods, Wehrfritz earlier proved in \cite[Thm. 5]{W} linearity of amalgamations of free groups along cyclic subgroups generated by primitive elements. Linearity of certain HNN extensions into $\mathsf{SL}(2,\mathbb{C})$ was also proved by Button in \cite[Thm. 6.1]{Button}.
\medskip

For our proof of Theorem \ref{t:LinearityBis} we need the following proposition.

\begin{proposition} \label{nonzero} Let $\Gamma$ be a torsion free word hyperbolic group and $\left \langle w \right \rangle$ be a maximal cyclic subgroup of $\Gamma$. Suppose that $\rho:\Gamma \rightarrow \mathsf{GL}(d, \mathbb{R})$ is a projective Anosov representation. Then there exists $h \in \mathsf{GL}(d, \mathbb{R})$ with the property: for every element $g \in \Gamma \smallsetminus \left \langle w \right \rangle$, the $(1,1), (1,d),(d,1)$ and $(d,d)$ entries of the matrix $h\rho(g)h^{-1}$ are non-zero. \end{proposition}

\begin{proof} Let $\xi_{\rho}:\partial_{\infty}\Gamma \rightarrow {\bf P}(\mathbb{R}^d)$ and $\xi^{-}_{\rho}:\partial_{\infty}\Gamma \rightarrow \mathsf{Gr}_{d-1}(\mathbb{R}^d)$ be the continuous $\rho$-equivariant Anosov limit maps of $\rho$. By the transversality of the Anosov limit maps we may find $h \in \mathsf{GL}(d,\mathbb{R})$ such that $\xi_{\rho}(w^{+})=h^{-1}\left \langle e_1 \right \rangle$, $\xi_{\rho}(w^{-})=h^{-1}\left \langle e_d \right \rangle$, $\xi^{-}_{\rho}(w^{+})=h^{-1}\left \langle e_d \right \rangle^{\perp}$ and $\xi^{-}_{\rho}(w^{-})=h^{-1}\left \langle e_1 \right \rangle^{\perp}.$ The matrices $\rho(w)$ and $\rho(w^{-1})$ are both 1-proximal and since $\xi_{\rho}$ and $\xi_{\rho}^{-}$ are dynamics preserving we may write $$\rho(w)=h^{-1}\begin{pmatrix}
s_{1} &0  &0 \\ 
 0& A &0 \\ 
0 & 0 & s_{d}
\end{pmatrix}h,$$ for some matrix $A \in \mathsf{GL}(d-2, \mathbb{R})$, such that $$|s_1|>\lambda_1(A) \geq \lambda_{d-2}(A)>|s_d|.$$ We remark that since $\Gamma$ is torsion-free and $\left \langle w \right \rangle$ is a maximal cyclic subgroup of $\Gamma$, the stabilizer of $w^{\pm}$ in $\Gamma$ (under the action of $\Gamma$ in $\partial_{\infty}\Gamma$) is the cyclic group $\langle w \rangle$ (see \cite{Gromov}). In particular, the cyclic group $\langle w \rangle$ is self-normalizing. Moreover, note that if $g \in \Gamma$ and $gw^{\pm}g^{-1} \in \{w^{+},w^{-}\}$, $gwg^{-1}$ has to be in the stabilizer of $w^{+}$ or $w^{-}$ under the action of $\Gamma$ on $\partial_{\infty}\Gamma$. Hence, $gwg^{-1}=w^{\pm 1}$, $g$ normalizes $\left \langle w \right \rangle$ and hence $g \in \langle w \rangle$. Therefore, we deduce that for every $g \in \Gamma \smallsetminus \left \langle w \right \rangle$, the intersection $\big \{gw^{+}, gw^{-}\} \cap \{w^{+},w^{-}\big \}$ is empty. \par  Finally, since the maps $\xi_{\rho}$ and $\xi^{-}_{\rho}$ are transverse, we have $$\rho(g)\xi_{\rho}(w^{+}) \oplus \xi^{-}_{\rho}(w^{\pm})=\rho(g)\xi_{\rho}(w^{-}) \oplus \xi^{-}_{\rho}(w^{\pm})=\mathbb{R}^d.$$ It follows that the $(1,1),(1,d),(d,1)$ and $(d,d)$ entries of $h\rho(g)h^{-1}$ are non-zero. \end{proof}

\noindent {\em Proof of Theorem \ref{t:LinearityBis}.} We split the proof of the theorem in two parts:

\noindent  {\em Construction of a family of representations $\big\{\pi_q:\Gamma \ast_{\langle w \rangle}\rightarrow G\big\}_{q>1}$.} By \cite[Lem. 3.18]{GW} we may assume that the inclusion of $\Gamma$ in $G$, $\rho:\Gamma \xhookrightarrow{} G$, is $\theta$-Anosov for some $\theta \subset \Delta$ which is stable under the opposition involution, i.e. $\theta^{\ast}=\theta$. Up to conjugating $\rho$ by an element $g \in G$, we may assume that the attracting fixed point of $\rho(w)$ (resp. $\rho(w)^{-1}$) in $G/P_{\theta}^{+}$ (resp. $G/P_{\theta}^{-}$) is the coset $P_{\theta}^{+}$ (resp. $P_{\theta}^{-}$) and $\rho(w) \in L_{\theta}$.\par By Proposition \ref{higherdimension}, there exists an irreducible $\theta$-proximal representation $\tau: G \rightarrow \mathsf{GL}(d, \mathbb{R})$ such that $\tau \circ \rho$ is a projective Anosov representation. Let $\chi_1,\ldots, \chi_{d}\in \mathfrak{a}^{\ast}$ be the restricted weights of $\tau$. Up to conjugating $\tau$, we may assume that $V^{\chi_i}=\langle e_i \rangle$ for $1 \leq i \leq d$. Observe that the dual representation $\tau^{\ast}$ is also $\theta$-proximal and its weights are $-\chi_1,\ldots,-\chi_d$. Up to conjugating $\tau$ by a permutation matrix of $\mathsf{O}(d)$, we may assume that the highest weight space of $\tau$ and $\tau^{\ast}$ are $V^{\chi_1}=\langle e_1 \rangle$ and $V^{\chi_d}=\langle e_d \rangle$ respectively.  

\par Let us now fix a vector $H_0 \in \overline{\mathfrak{a}}^{+} \cap \mathfrak{a}_{\theta}$ such that $$\alpha(H_0) \in \mathbb{N}^{\ast}, \forall \ \alpha \in \theta \ \ \textup{and} \ \ \alpha(H_0)=0, \ \forall \ \alpha \in \Delta \smallsetminus \theta.$$ Let $q>1$. The matrix $\tau \big(  \exp(\log(q) H_0)\big)=\exp\big(\log (q) d\tau (H_0)\big)$ has the following properties:
\begin{itemize}
\item[(i)] $\tau \big(\log (q) \exp(H_0)\big)$ is the diagonal matrix $\textup{diag}\big(q^{\chi_1(H_0)},\ldots,q^{\chi_{d}(H_0)}\big)$.
\medskip

\item[(ii)] $\chi_1(H_0)-\chi_i(H_0)\in \mathbb{N}^{\ast}$ and $\chi_d(H_0)-\chi_j(H_0)\in \mathbb{N}^{\ast}$ for every $2 \leq i \leq d$  and $1 \leq j \leq d-1$. The attracting fixed points of $\tau\big( \exp( \log(q) H_0)\big)$ and $\tau \big(\exp(-\log(q) H_0)\big)$ in ${\bf P}(\mathbb{R}^d)$ are the lines $[e_1]$ and $[e_d]$ respectively.
\medskip

\item[(iii)] $\exp \big( \pm \log(q) H_0 \big)$ commutes with $\rho(w)\in L_{\theta}$. \end{itemize} 
\medskip

By Proposition \ref{univpropHNN} we conclude that there exists a well defined group homomorphism ${\pi_q:\Gamma \ast_{\left \langle w \right \rangle}\rightarrow G}$ defined as follows: \begin{align*} \pi_q(t)&=\exp \big(\log(q)H_0 \big)\\ \pi_q(\gamma)&=\rho(\gamma), \ \gamma \in \Gamma. \end{align*} 

\medskip

\noindent {\em Injectivity of $\pi_q:\Gamma \ast_{\langle w \rangle}\rightarrow G$ for generic values of $q>1$.}  
We recall that the attracting fixed point of $\rho(w)$ in $G/P_{\theta}^{+}$ is the coset $P_{\theta}^{+}$ and its repelling fixed point in $G/P_{\theta}^{-}$ is $P_{\theta}^{-}$. It follows by Remark \ref{Plucker} that the attracting fixed point of $\tau(\rho(w))$ $\big($resp. $\tau^{\ast}(\rho(w))$$\big)$ in ${\bf P}(\mathbb{R}^d)$ is the line $\iota_{\tau}^{+}(P_{\theta}^{+})=\big[V^{\chi_{\tau}}\big]=[e_1]$ $\big($resp. $\iota_{\tau^{\ast}}^{+}(P_{\theta}^{+})=\big[V^{\chi_{\tau^{\ast}}}\big]=[e_d]$$\big)$. The repelling fixed point of $\tau^{\ast}(\rho(w))$ (resp. $\tau(\rho(w))$) in $\mathsf{Gr}_{d-1}(\mathbb{R}^d)$ is the $(d-1)$-plane $\iota_{\tau^{\ast}}^{-}(P_{\theta}^{-})=V^{<\chi_{\tau^{\ast}}}=\langle e_d \rangle ^{\perp}$ (resp. $\iota_{\tau}^{-}(P_{\theta}^{-})=V^{< \chi_{\tau}}=\langle e_1 \rangle^{\perp}$). We deduce from Proposition \ref{nonzero} and its proof that for every $\gamma \in \Gamma \smallsetminus \langle w\rangle$ the $(1,1),(1,d),(d,1)$ and $(d,d)$ entries of the matrix $\tau(\rho(\gamma))$ are non-zero. 

\par Let $\mathbb{F}$ be the finitely generated subfield of $\mathbb{R}$ spanned by the entries of the elements of $\tau(\rho(\Gamma))$. Let us chose $q>1$ to be transcendental over the field $\mathbb{F}$. We claim that the representation \hbox{$\pi_q: \Gamma \ast_{\langle w \rangle}\rightarrow G$} is faithful. Suppose that $h \in \Gamma \ast_{\langle w \rangle}$ is a non-trivial element. If $h$ lies in a conjugate of $\Gamma$ we obviously have $\tau(\pi_q(h)) \neq \textup{I}_d$ since $\Gamma$ is torsion free and $\tau \circ \rho: \Gamma \rightarrow \mathsf{GL}(d,\mathbb{R})$ is projective Anosov and faithful. If $h \in \Gamma \ast_{\langle w \rangle}$ does not lie in a conjugate of $\Gamma$, by Britton's lemma for HNN extensions (see Section \ref{ss:Amalgamated products and HNN}), $h$ is conjugate to a product of the form $$h_{k}=t^{p_1}g_1t^{p_2}g_2 \cdot \cdot \cdot t^{p_{k}}g_{k},$$ where $g_{j}\in \Gamma \smallsetminus \left \langle w \right \rangle$ for $1 \leq j \leq k$ and $p_1 \neq 0$. We may assume that $p_1>0$ and we will show that $\tau(\pi_q(h_k))$ is not a scalar multiple of $\mathrm I_d$. If $p_1<0$ similar arguments will apply. We may write $$\tau(\pi_q(t))^{p_i}=\textup{diag}\Big( q^{p_i \chi_1(H_0)},\ldots,q^{p_i \chi_d(H_0))}\Big)=q^{m_i}A_{i}$$ where 
$$A_{i} \equaldef \left\{\begin{matrix}
\textup{diag}\Big(q^{p_i(\chi_1(H_0)-\chi_d(H_0))},\ldots,1\Big) & \textup{if} \ p_i>0 \\
\textup{diag}\Big(1,\ldots,q^{p_i(\chi_d(H_0)-\chi_1(H_0))}\Big) & \textup{if} \ p_i<0.
\end{matrix}\right.$$ and $m_i=p_i \chi_d(H_0)$ if $p_i>0$ and $m_i=p_i \chi_1(H_0)$ if $p_i<0$.
In particular, since $\pi_q(\gamma)=\rho(\gamma)$ for every $\gamma \in \Gamma$, we may write \begin{align*} \tau(\pi_q(h_k))&=\tau(\pi_q(t))^{p_1}\tau(\rho(g_1))\cdot \cdot \cdot \tau(\pi_q(t))^{p_k}\tau(\rho(g_k))\\
&=q^{s} \big(A_1 \tau(\rho(g_1))\big)\cdot \cdot \cdot \big(A_k \tau(\rho(g_k))\big),\end{align*} where $s=\sum_{i=1}^{k}m_i$. Let us also set $s_i \equaldef p_{i}\big(\chi_1(H_0)-\chi_d(H_0)\big)$ for $1 \leq i \leq k$. We have already seen that for every $1 \leq j \leq k$ the $(1,1),(1,d),(d,1)$ and $(d,d)$ entries of $\tau(\rho(g_j))$ are non-zero. We first check that $A_1 \tau(\pi_q(g_1))$ is a matrix whose $(1,d)$ and $(d,1)$ entries are polynomials in $q$ of degree $s_1$, for $2 \leq i \leq d-1$ the $(1,i)$ entry is a polynomial of degree at most $s_1$ and the remaining entries are polynomials of degree at most $s_1-1$.

Next we mutliply with the matrix $A_2 \tau(\pi(g_2))$. There are two cases to consider:
\medskip

\noindent {\em Case 1}. Suppose that $p_2>0$ (and $s_2>0$). The $(1,1)$ and $(1,d)$ entries of $A_2 \tau(\rho(g_2))$ is a polynomial in $q$ of degree $s_2$, the remaining $(1,i)$ entries have degree at most $s_2$ and all the other entries have degree at most $s_2-1$. We see that the $(1,1)$ (resp. $(1,d)$) entry of $A_1\tau(\rho(g_1))A_2\tau(\rho(g_2))$ is obtained by multiplying the $(1,1)$ entry of $A_1 \tau(\rho(g_1))$ with the $(1,1)$ (resp. $(1,d)$) entry of $A_2 \tau(\rho(g_2))$ plus we add some terms of degree at most $s_1+s_2-1$. With this observation, we see that the $(1,1)$ and $(1,d)$ entries of $A_1 \tau(\rho(g_1))A_2\tau(\rho(g_2))$ are polynomials in $q$ of degree $s_{1}+s_{2}$. For $2 \leq i \leq d-1$, the entries $(1,i)$ of $A_1 \tau(\rho(g_1))A_2\tau(\rho(g_2))$ have degree at most $s_1+s_2$ and the remaining entries have degree at most $s_1+s_2-1$.
\medskip

\noindent {\em Case 2}. Suppose that $p_2<0$ (and $s_2<0$). The product $A_2 \tau(\rho(g_2))$ has its $(d,1)$ and $(d,d)$ entry as a polynomial in $q$ of degree $|s_{2}|$, all the other entries $(d,i)$ are of degree at most $|s_{2}|$ and the remaining entries have degree at most $|s_2|-1$. We check that the $(1,1)$ (resp. $(1,d)$) entry of $A_1\tau(\rho(g_1))A_2\tau(\rho(g_2))$ is obtained by multiplying the $(1,d)$ entry of $A_1 \tau(\rho(g_1))$ with the $(d,1)$ (resp. $(d,d)$) entry of $A_2 \tau(\rho(g_2))$ plus we add some terms of degree at most $s_1+|s_2|-1$.  In this case, we deduce that the $(1,1)$ and $(1,d)$ entries of $A_1 \tau(\rho(g_1))A_2\tau(\rho(g_2))$ are polynomials in $q$ of degree $s_{1}+|s_{2}|$. The remaining entries $(1,i)$ for $2 \leq i \leq d-1$ are of degree at most $s_{1}+|s_{2}|$ and all other entries are polynomials of degree at most $s_1+|s_{2}|-1$.

By induction, one shows that the $(1,1)$ and $(1,d)$ entries of the product $A_1 \tau(\rho(g_1)) \cdot \cdot \cdot A_k \tau(\rho(g_k))$ are polynomials  in $\mathbb{F}[q]$ of degree $$d_{k}=\sum_{i=1}^{k} |s_i|=\big(\chi_1(H_0)-\chi_d(H_0)\big)\sum_{i=1}^{k} |p_i|,$$ for every $2 \leq i \leq d-1$ the $(1,i)$ entry is a polynomial of degree at most at most $d_k$ and all the other entries have degree at most $d_k-1$ in $q$. Since $q$ was chosen to be transcendental over $\mathbb{F}$, It follows that $\tau(\pi(h))$ is not a scalar mutliple of $\mathrm I_d$. \par Finally, we conclude that $\tau \circ \pi_q:\Gamma \ast_{\langle w \rangle} \rightarrow \mathsf{GL}(d,\mathbb{R})$ (and hence $\pi_q$) is a faithful representation.  \hfill $\qed$ \\

As a corollary of the method of proof of Theorem \ref{t:LinearityBis} we have:

\begin{corollary} Let $G$ be a linear semisimple Lie group. Fix $\theta \subset \Delta$ a subset of simple restricted roots of $G$ and let $\Gamma_1$ and $\Gamma_2$ be two torsion-free word hyperbolic groups. Let $\langle w_1\rangle$ and $\langle w_2 \rangle$ be two maximal cyclic subgroups of $\Gamma_1$ and $\Gamma_2$ respectively. Suppose that $\rho_1:\Gamma_1 \rightarrow G$ and $\rho_2:\Gamma \rightarrow G$ are $\theta$-Anosov representations and $\rho_1(w_1)=\rho_2(w_2)$. Then there exists $h \in G$ such that the subgroup $\big \langle \rho_1(\Gamma_1), h\rho_2(\Gamma)h^{-1}\big \rangle$ of $G$ is isomorphic to the amalgamated product $\Gamma_1 \ast_{w_1=w_2}\Gamma_2$.  \end{corollary}

\begin{proof} We keep the notation from the previous proof. Up to conjugating both $\rho_1$ and $\rho_2$  by some element $g \in G$, we may assume that $\rho_1(w_1)=\rho_2(w_2) \in L_{\theta}$. As previously, let $\tau:G \rightarrow \mathsf{GL}(d,\mathbb{R})$ be an irreducible $\theta$-proximal representation such that $\tau \circ \rho_i:\Gamma_i \rightarrow \mathsf{GL}(d,\mathbb{R})$ is a projective Anosov representation for $i=1,2$. We choose $q>1$  transcendental over the finitely genrated field $\mathbb{F}'$ generated by the entries of the matrices in $\tau(\rho_1(\Gamma)) \cup \tau(\rho_2(\Gamma))$. We may assume that the matrix $h_{q} \equaldef \exp \big(\log(q)d\tau( H_0))\big) $ is a diagonal matrix of the form $\textup{diag}\big(q^{\chi_1(H_0)},\ldots,q^{\chi_d(H_0)}\big)$, where $\chi_1(H_0)-\chi_{i}(H_0) \in \mathbb{N}^{\ast}$ for every $2 \leq i \leq d$. Then there is a well defined group homomorphism $\pi_{q}':\Gamma_1 \ast_{w_1=w_2}\Gamma_2 \rightarrow G$ such that \begin{align*} \pi_{q}'(\gamma_1)& =\rho_1(\gamma_1), \ \gamma_1 \in \Gamma_1 \\  \pi_{q}'(\gamma_2)&=h_q  \rho_2(\gamma_2) h_{q}^{-1}, \ \gamma_2 \in \Gamma_2. \end{align*}  By Proposition \ref{nonzero}, the $(1,1),(1,d),(d,1)$ and $(d,d)$ entries of $\rho_i(\gamma_i)$ are non-zero for $\gamma_i \in \Gamma_i \smallsetminus \langle w_i \rangle$ and $i=1,2$. We similarly check that for a word $h_{k}=g_{11}g_{21} \cdots g_{1k}g_{2k}$, $g_{ij}\in \Gamma_{i}\smallsetminus \langle w_i \rangle$, which is not in a conjugate of $\Gamma_{1}$ or $\Gamma_2$, the $(1,1)$ and $(1,d)$ entries of $\tau(\pi_{q}'(h_k))$ are of the form $q^{r_k}f(q)$, where $r_k \in \mathbb{R}$ and $f(q) \in \mathbb{F}[q]$ is a polynomial of degree $k(\chi_1(H_0)-\chi_d(H_0))$. In particular, $\pi_{q}'$ is injective and the group $\big \langle \rho_1(\Gamma_1), h_{q} \rho_2(\Gamma_2) h_{q}^{-1} \big \rangle$ is isomorphic to $\Gamma_1 \ast_{w_1=w_2}\Gamma_2$.
\end{proof}

\section{Superrigidity and arithmeticity} \label{s:superrigidity}

The renowned superrigidity theorem of Margulis \cite{Mar} states that linear representations of an irreducible lattice in a real semisimple linear group $G$ of rank at least $2$ essentially extend to the whole group $G$. The theorem was extended by Corlette \cite{Corlette} to representations of lattices in the quaternionic groups $\Sp(k,1)$, $k\geq 2$ and the exceptional group $\mathrm F_4^{(-20)}$ (of rank $1$).

Let $G$ be a real simple linear group which is either isogeneous to $\Sp(k,1)$ or $\textup{F}_4^{(-20)}$ or of rank at least $2$, and let $\Gamma$ be a lattice in $G$. We first state a geometric version of these superrigidity theorems.

\begin{theorem}  \label{t:SuperRigidityGeodesic} 
Let $H$ be a real semisimple linear group and $\rho$ a representation of $\Gamma$ into $H$. Then there exists a $\rho$-equivariant map
\[f:G/K \to H/L\]
which is  totally geodesic. \textup{(}Here $G/K$ and $H/L$ denote respectively the symmetric spaces of $G$ and $H$.\textup{)}
\end{theorem}

\begin{remark} \normalfont{Since $G$ is simple, the symmetric space $G/K$ is irreducible. Hence the totally geodesic map $f$ is either constant \textup{(}in which case $\rho$ takes values in a compact subgroup of $H$\textup{)} or an embedding which is isometric up to scaling one of the symmetric metrics.} \end{remark}

\medskip

Margulis (in higher rank) and Gromov--Schoen \cite{Gromov-Schoen} (for $G=\Sp(k,1)$ or $\mathrm F_4^{(-20)}$) also proved Theorem \ref{t:SuperRigidityGeodesic}  for linear representations of $\Gamma$ over (complete) valued fields. There, the Bruhat--Tits building of $H$ plays the role of the symmetric space, and totally geodesic maps from a Riemannian symmetric space are constant, implying that every representation of $\Gamma$ over such a field has bounded image.

This stronger form of superrigidity has many consequences. It implies in particular that all lattices under consideration are \emph{arithmetic}:

\begin{theorem} \label{t:Arithmeticity}
Given $G$ and $\Gamma$ as above, there exists a semisimple linear algebraic group $\mathrm G$ over $\mathbb Q$ and a smooth surjective morphism $\phi: \mathrm G(\Real) \to G$ with compact kernel such that $\Gamma$ is comensurable to $\phi(\mathrm G(\mathbb Z))$.
\end{theorem}

Finally, these superrigidity results essentially classify the linear representations of lattices in higher rank. Let us now assume (without loss of generality by Theorem \ref{t:Arithmeticity}) that $\Gamma$ is commensurable to $\mathrm G(\mathbb Z)$, where $\mathrm G$ is a semisimple algebraic group over $\mathbb Q$ such that $\mathrm G(\mathbb R)$ is isogeneous to a product of $G$ with a compact group.

\begin{theorem} \textup{(}\cite[Cor. 16.4.1]{Morris}\textup{)} \label{t:RepresentationsSuperRigidLattices}
Let $\rho: \Gamma \to \GL(n,\mathbb R)$ be a linear representation. Then there exists a smooth morphism $\phi: \mathrm G(\Real) \to \GL(n,\mathbb R)$ which coincides with $\rho$ on a finite index subgroup of $\Gamma$.
\end{theorem}

\begin{remark} \label{rem:Representations finite characteristic}
\normalfont{Let us finally note that superrigidity over valued fields also implies that linear representations of $\Gamma$ into $\GL(n,{\bf k})$ have finite image when ${\bf k}$ is a field of positive characteristic \textup{(}see the proof in \cite[Thm. 8.1]{Kap}\textup{)}}.
\end{remark}

\section{Non-linear and indiscrete groups} \label{s:proofs}

\subsection{Non-linear amalgamations}

Here we prove the three non-linearity results stated in the introduction. As in previous examples, see \cite{Kap,CST}, superrigid lattices will be the starting points of our constructions of non-linear groups. We first prove Theorem \ref{t:NonLinearAmalgamatedCyclicIntro} which we recall here.

\medskip

\noindent {\bf Theorem \ref{t:NonLinearAmalgamatedCyclicIntro}} {\em Let $\Gamma_1$ and $\Gamma_2$ be lattices in $\Sp(k,1)$, $k \geq 2$, and $\langle w_i \rangle$ be a cyclic subgroup of $\Gamma_i$ for $i=1,2$. Assume that $w_1 \in \Gamma_1$ and $w_2 \in \Gamma_2$ have different translation lengths in the symmetric space of $\Sp(k,1)$. Then every linear representation of $\Gamma_1\ast_{w_1=w_2} \Gamma_2$ restricted on $\Gamma_1$ and $\Gamma_2$ has finite image.}
\medskip

For an isometry $g\in \mathsf{Sp}(k,1)$ we denote by $\ell_{\mathbb{H}{\bf H}^k}(g)$ the translation length for its action on $\mathbb{H}{\bf H}^k$. Let us note that we equip the quaternionic hyperbolic space $\mathbb{H}{\bf H}^k$ with the negatively curved Riemannian metric induced by a scalar multiple of the Killing form (on the symmetric part of the Lie algebra of $\mathsf{Sp}(k,1)$) such that the hyperbolic isometry \[g(t) \equaldef \left( \begin{matrix} \textup{I}_{k-1} & 0 & 0\\
                                                                        0 & \cosh(t) & \sinh(t)\\
								0 &  \sinh(t) & \cosh(t)
					                                \end{matrix} \right)~\]
satisfies $\ell_{\mathbb{H}{\bf H}^k}(g(t)) = t$ or every $t\geq 0$.

For the proof we will need the following lemma. Recall that for an element $g\in \mathsf{GL}(d,\mathbb{C})$, $\lambda_1(g)\geq \cdots \geq \lambda_d(g)$ are the moduli of the eigenvalues of $g$ in decresing order and $\overrightarrow{\ell} =(\log \lambda_1,\ldots, \log \lambda_r): \GL(r,\Real) \to \Real^r$ denotes the Jordan projection. A matrix $g$ is called $1$-proximal if $\lambda_1(g)>\lambda_2(g)$.

\begin{lemma} \label{l:JordanProjectionsTranslationLength}
Let $k \geq 2$ and $\rho: \mathsf{Sp}(k,1) \rightarrow \GL(r,\CC)$ be a non-trivial continuous representation. Let $1 \leq i\leq r-1$ be the largest index such that $\log \lambda_i (\rho(g)) =\log \lambda_1 (\rho(g))$ for every $g \in \mathsf{Sp}(k,1)$. Then for all $g\in \mathsf{Sp}(k,1)$ we have
\[ \log \frac{\lambda_i(\rho(g))}{\lambda_{i+1}(\rho(g))} = \ell_{\mathbb{H}{\bf H}^k}(g).\]\end{lemma}

\begin{proof} Note that by the definition of the index $i\in \mathbb{N}$ the representation $\wedge^i \rho:\mathsf{Sp}(k,1)\rightarrow \mathsf{GL}(\wedge^i \mathbb{C}^r)$ is $1$-proximal.
Since $\Sp(k,1)$ has real rank $1$ it is enough to determine the eigenvalues of the matrix $\wedge^i \rho(g(t))$.
Now note that the restriction $\wedge^i \rho: \SU(2,1)\times \{ \textup{I}_{k-2}\}\rightarrow \mathsf{GL}(\wedge^i \mathbb{C}^r)$ can be extended to a complex semisimple representation $\psi:\SL(3,\CC)\rightarrow \mathsf{GL}(\wedge^i \mathbb{C}^r)$ which decomposes as a direct product $$\psi=\psi_1 \times \cdots \times \psi_{p}$$ where $\big\{\psi_i:\mathsf{SL}(3,\mathbb{C})\rightarrow \mathsf{SL}(V_i)\big\}_{i=1}^{p}$ are irreducible complex representations and $\wedge^i\mathbb{C}^d=V_1\oplus \cdots \oplus V_r$. Note that $\psi(g(t))$ is $1$-proximal and its attracting fixed point in $\mathbb{P}(\mathbb{C}^d)$ necessarily lies in $\bigcup_{i=1}^{r}\mathbb{P}(V_i)$, say in $\mathbb{P}(V_1)$. In particular, $\psi_{1}$ is $1$-proximal and $\log \lambda_1(\psi_1(g(t))>\log \lambda_1(\psi_j(g(t))$ for $2\leq j \leq p$ and $t>0$.
By using the representation theory of $\mathsf{SL}(3,\mathbb{C})$ or \cite[Lem. 3.7]{GGKW}, one verifies that $\log \lambda_1(\psi_j(g(t)))$ is an integral multiple of $t$ for every $1\leq j \leq p$ and also $\log \frac{\lambda_1(\psi_1g(t)))}{ \lambda_2(\psi_1(g(t)))}= t$.\footnote{Here, it is important that $k\geq 2$ in order to guarantee that $1$ is indeed an eigenvalue of $g(t)$.} Moreover, observe that $$\log \frac{\lambda_1(\psi(g(t)))}{\lambda_2(\psi(g(t)))}=\min\Big\{ \log \frac{\lambda_1(\psi_1(g(t)))}{\lambda_2(\psi_1(g(t)))}, \log \frac{\lambda_1(\psi_1(g(t)))}{\lambda_1(\psi_2(g(t)))}, \ldots, \log \frac{\lambda_1(\psi_1(g(t)))}{\lambda_1(\psi_p(g(t)))} \Big\}$$ and therefore it follows that $\log \frac{\lambda_i(\rho(g(t)))}{\lambda_{i+1}(\rho(g(t)))}=\log \frac{\lambda_1(\psi(g(t)))}{ \lambda_2(\psi(g(t)))}=t$. The conclusion follows.\end{proof}
\medskip

\noindent {\em Proof of Theorem \ref{t:NonLinearAmalgamatedCyclicIntro}.}
Let $L$ be a field and $\rho: \Gamma_1 \ast_{g_1=g_2} \Gamma_2 \to \GL(d,L)$ be a linear representation. Assume first that $L$ has positive characteristic. Then $\rho|_{\Gamma_1}$ has necessarily finite image, see Remark \ref{rem:Representations finite characteristic} and \cite{Kap}. Hence $\rho$ cannot be injective.

\par We now assume that $L$ has characteristic $0$. Since $\Gamma_1 \ast_{g_1=g_2} \Gamma_2$ is finitely generated, we may assume without loss of generality that $L$ is finitely generated over $\mathbb Q$. If $\rho|_{\Gamma_1}$ has infinite image, then there exists a representation $\tau: \GL(d,L) \to \GL(r,\mathbb R)$ such that $\tau \circ \rho|_{\Gamma_1}$ has unbounded image (see \cite[Thm. 3.1]{CST}). By Corlette's superrigidity theorem (see Theorem \ref{t:SuperRigidityGeodesic}), there exists a continuous non-trivial representation $\rho_1: \Sp(k,1)\to \GL(r,\Real)$ such that \[\overrightarrow{\ell}(\rho_1(g)) =\overrightarrow{\ell}(\rho(g))\]
for every $g\in \Gamma_1$. Up to taking an exterior power of $\rho$ and $\rho_1$, we may also assume that $\rho_1$ (and hence $\rho$) is $1$-proximal. In particular, $\lambda_1(\rho_1(g_1))>\lambda_2(\rho_1(g_1))$.

Now we observe that $\overrightarrow{\ell}(\rho(g_1))=\overrightarrow{\ell}(\rho(g_2)) \neq 0$ and hence the image of the restriction $\rho|_{\Gamma_2}$ is also unbounded and contains a $1$-proximal elemment. By Corlette's superrigidity there exists a continuous proximal representation $\rho_2: \Sp(k,1) \to \GL(r,\Real)$ such that \[\overrightarrow{\ell}(\rho_2(g)) =\overrightarrow{\ell}(\rho(g))\]
for every $g\in \Gamma_2$. Finally, we have $$\log \frac{\lambda_1(\rho_1(g_1))}{\lambda_2(\rho_1(g_1))}= \log \frac{\lambda_1(\rho(g_1))}{\lambda_2(\rho(g_1))}=\log \frac{\lambda_1(\rho(g_2))}{\lambda_2(\rho(g_2))}=\log \frac{\lambda_1(\rho_2(g_2))}{\lambda_2(\rho_2(g_2))}.$$
By using Lemma \ref{l:JordanProjectionsTranslationLength} for the $1$-proximal representations $\rho_1$ and $\rho_2$ of $\mathsf{Sp}(k,1)$ and $i=1$, we obtain that the translation lengths of $g_1$ and $g_2$ on $\mathbb{H}{\bf H}^k$ have to be equal. However, this contradicts the hypothesis that $$\ell_{\mathbb{H}{\bf H}^k}(g_1)\neq \ell_{\mathbb{H}{\bf H}^k}(g_2).$$ Finally, we conclude that every linear representation of the amalgamated product $\Gamma_1\ast_{g_1=g_2}\Gamma_2$ over a field $L$ restricted on either $\Gamma_1$ or $\Gamma_2$ has finite image. $\qed$
\medskip

Our remaining non-linearity results will follow from a general lemma about representations of superrigid lattices. As in Section \ref{s:superrigidity}, let $G$ be a real semisimple linear group which is either of rank at least two or isogeneous to $\Sp(k,1), k\geq 2$, or $\mathrm F_4^{(-20)}$ and let $\Gamma$ be a lattice in $G$.
\medskip

\begin{lemma} \label{l:RepresentationsCoincidingOnW}
Let $W$ be a subgroup of $G$. Let $\rho_1: \Gamma \rightarrow \GL(r,{\bf k})$ and $\rho_2:\Gamma \to \GL(r,{\bf k})$ be two linear representations of $\Gamma$ over a field ${\bf k}$ which coincide on $W$. Then $\rho_1$ and $\rho_2$ coincide on a finite index subgroup of $\Gamma \cap \overline{W}^Z$. \textup{(}Here, $\overline{W}^Z$ denotes the Zariski closure of $W$ in $G$.\textup{)}
\end{lemma}

\begin{proof}
Assume first that ${\bf k}$ has positive characteristic. Then $\rho_1$ and $\rho_2$ have finite image by Remark \ref{rem:Representations finite characteristic}, hence they are both trivial on a finite index subgroup of $\Gamma$.

Let us now assume that ${\bf k}$ has characteristic $0$. Since $\Gamma$ is finitely generated, $\rho_1$ and $\rho_2$ have their image in a finitely generated extension of $\mathbb Q$ (the extension generated by all the coefficients of the image of a finite generating subset of $\Gamma$). We can thus assume that ${\bf k}$ embeds in $\mathbb C$, and after composing with the restriction of scalars $ \GL(r,\mathbb C) \to \GL(2r,\mathbb R)$ we can thus restrict to the case where ${\bf k}= \mathbb R$.

By Theorem \ref{t:Arithmeticity}, there exists a semisimple linear algebraic group $\mathrm G$ over $\mathbb Q$ and a smooth morphism $\phi: \mathrm{G}(\mathbb R) \to G$ with compact kernel such that $\Gamma$ is commensurable to $\phi(\mathrm G(\mathbb Z))$. Up to passing to a finite index subgroup and a finite cover, we can assume that $\mathrm{G}$ is algebraically connected and simply connected (equivalently, that $\mathrm{G}(\mathbb C)$ is connected and simply connected). This implies that the morphism $\phi$ is algebraic (see \cite{Morris}). If we replace $\Gamma$ by $\phi^{-1}(\Gamma)$, $W$ by $\phi^{-1}(W)$ and $\rho_1$ and $\rho_2$ by $\rho_1\circ \phi$ and $\rho_2 \circ \phi$, we are thus reduced to the case where $G= \mathrm{G}(\mathbb R)$ and $\Gamma$ is commensurable to $\mathrm G(\mathbb Z)$.

Now, by Theorem \ref{t:RepresentationsSuperRigidLattices}, there exist $\psi_1$ and $\psi_2: \mathrm{G}(\mathbb R) \to \GL(r,\mathbb R)$ that coincide with $\rho_1$ and $\rho_2$ on a finite index subgroup of $\Gamma$. Since $\mathrm{G}$ is algebraically simply connected, $\psi_1$ and $\psi_2$ are algebraic over $\mathbb R$. Now, since $\rho_1$ and $\rho_2$ coincide on $W$, $\psi_1$ and $\psi_2$ coincide on a finite index subgroup of $W$. Since they are algebraic morphisms, they coincide on a finite index subgroup of $\overline{W}^Z$, and $\rho_1$ and $\rho_2$ coincide on a finite index subgroup of $\Gamma \cap \overline{W}^Z$.\end{proof}

We obtain the following corollary which immediately implies Theorem \ref{t:NonLinearDoubleIntro}.

\begin{corollary}  \label{c:NonLinearDoublesArithmeticLattices}
If $W$ is not a lattice in $\overline{W}^Z$, then $\Gamma \ast_W \Gamma$ is not linear.
\end{corollary}

\begin{proof} Let us first note that $\Gamma \cap \overline{W}^Z$ is a lattice in $\overline{W}^Z$. Indeed, we can restrict to the case where $G= \mathrm G(\mathbb R)$ and $\Gamma$ is commensurable to $\mathrm G(\mathbb Z)$. Then $\overline{W}^Z$ is defined over $\mathbb Q$. Since $W \subset \mathrm G(\mathbb Z)$, every $\mathbb Q$-character of $\overline{W}^Z$ is virtually trivial on $W$, hence trivial on $\overline{W}^Z$. Therefore, $\Gamma \cap \overline{W}^Z$ is a lattice in $\overline{W}^Z$ by Borel--Harish-Chandra's theorem (see for instance \cite[Th\'eor\`eme 5.4]{Benoist}).

Let $\rho:\Gamma \ast_W \Gamma \rightarrow \mathsf{GL}(d,{\bf k})$ be a linear representation of $\Gamma \ast_W \Gamma$. By Lemma \ref{l:RepresentationsCoincidingOnW} and the universal property of amalgamated products, $\rho$ factors through $\Gamma \ast_{W'} \Gamma$, where $W' \supset W$ is a finite index subgroup of $\Gamma \cap \overline{W}^Z$. Since $W$ is not a lattice in $\overline{W}^Z$, the group $W'$ strictly contains $W$ and hence Fact \ref{nonfaithful} implies that the representation $\rho$ is not faithful.
\end{proof}

\medskip

Theorem \ref{t:NonLinearitySLZ}, on the other side, follows from Corollary \ref{c:NonLinearDoublesArithmeticLattices} applied to a maximal cyclic subgroup $\langle w \rangle$ of $\SL(d,\mathbb Z)$ whose Zariski closure contains a real split torus of $\mathsf{SL}(d,\mathbb{R})$ of rank at least two. Such a subgroup exists by the following lemma, which can be found in \cite{PR}.

\begin{lemma} \label{l:ZariskiColsureCyclicSL(d,Z)}
There exists $g\in \SL(d,\Z)$ such that $\overline{\langle g \rangle}^Z$ is a real split maximal torus.
\end{lemma}

We prove this lemma here for completeness. Let us fix $K$ a totally real extension of $\mathbb Q$ of degree $d$ and $L$ its smallest Galois extension. Let $M$ be the Galois group of $L$ and $N$ the subgroup fixing $K$. Consider the algebraic group $\mathrm H = \mathrm{Res}_{K/\mathbb Q}(\mathrm G_m(K))$, where $\mathrm G_m(K)$ is the multiplicative group of $K$.

Let us fix an embedding of $L$ into $\mathbb R$. Then the group morphism
\begin{eqnarray*}
K^\times &\to & {\Real^\times}^{M/N}\\
a  & \mapsto & (\sigma(a))_{\sigma \in M/N}
\end{eqnarray*}
Induces an algebraic isomorphism from $\mathrm H(\mathbb R)$ to $\mathrm G_m(\mathbb R)^{M/N}$. The image of $K^\times$ is identified with the subgroup $\mathrm H(\mathbb Q)$. The group $M$ acts on $\mathrm G_m(\mathbb R)^{M/N}$ by permutation of the coordinates.

\begin{proposition} \label{p:Qtorus Galois invariant}
Let $\mathrm H'$ be a $\mathbb Q$-algebraic subgroup of $\mathrm H$. Then $\mathrm H'(\mathbb R) \subset \mathrm G_m(\mathbb R)^{M/N}$ is invariant under the action of $M$.
\end{proposition}

\begin{proof}
Since $\mathrm H'$ is an algebraic subgroup of a real split torus, we have
\[\mathrm H' = \bigcap_{\chi_{\vert \mathrm H'}  =1} \ker \chi~,\]
where the intersection is taken over all real characters of $\mathrm H$ that are trivial on $\mathrm H'$.

Now, let $\chi$ be a real character of $\mathrm H$ which is trivial on $\mathrm H'$. This character is given by
\[\chi(x_\sigma)= \prod_{\sigma \in M/N} x_\sigma^{m_\sigma}\]
For some $(m_\sigma) \in \Z^{M/N}$.
If $x$ belongs to $\mathrm H'(\mathbb Q)$, then $x=(\sigma(a))$ for some $a\in K^\times$, and for all $\eta \in M$ we have
\begin{eqnarray*}
\chi (\eta\cdot x) &=& \chi (\eta \sigma(a))\\
&=& \prod_{\sigma \in M/N} (\eta \sigma(a))^{m_\sigma}\\
&=& \eta (\chi(a)) = 1~.
\end{eqnarray*}
We conclude that $M \cdot \mathrm H'(\mathbb Q) \subset \mathrm H'(\mathbb R)$, hence $\mathrm H'(\mathbb R)$ is $M$-invariant.
\end{proof}

Now, let $O_K$ denote the ring of integers of $K$ and $\Lambda$ be the multiplicative group of integers of norm $1$ (which has index $2$ in $O_K^\times$). By Dirichlet's unit theorem, $\Lambda$ is a lattice in the group $\mathrm H_0(\mathbb R) = \{(a_\sigma)\in \mathrm G_m(\mathbb R)^{M/N} \mid \prod_{\sigma \in M/N} a_\sigma = 1\}$.

\begin{proposition} \label{l:ZariskiDenseCyclicSubgroupTorus}
There exists $a\in \Lambda$ such that $\langle a \rangle$ is Zariski dense in $\mathrm H'$.
\end{proposition}

\begin{proof}
Define
\[\begin{array}{cccc}
\mathrm e: & \Real^{M/N} & \to & \mathrm G_m(\Real)^{M/N} \\
& (x_\sigma)_{\sigma\in M/N} & \mapsto & (e^{x_\sigma})_{\sigma \in M/N}~.
\end{array}\]
The morphism $\mathrm e$ commutes with the action of $M$ on coordinates and induces a smooth isomorphism from the hyperplane $V= \{ \sum_\sigma x_\sigma = 0\}$ to the identity component of $\mathrm H_0(\mathbb R)$, so that the preimage of $\Lambda$ by $\mathrm e$ is a lattice in the vector space $V$.

The set $U$ of vectors $x\in V$ such that $\mathrm{Span} \{ \sigma\cdot x\}_{\sigma\in M} =V$ is a Zariski-open cone, which is non-empty since the vector $(d-1, -1,\ldots, -1)$ belongs to $U$. It must therefore contain an element $x$ (and in fact, most elements) of $\mathrm e^{-1}(\Lambda)$.

Let $H'$ be the Zariski closure of $\langle \mathrm e(x)\rangle$.
 By Proposition \ref{p:Qtorus Galois invariant}, $\mathrm{H}'(\mathbb R)$ contains all the $\mathrm e (\sigma \cdot x)$ for $\sigma \in M$. Let $\chi$ be a real character of $\mathrm H$ which is trivial over $\mathrm H'$. Then $\log \chi \circ \mathrm e$ is a linear form on $\Real^{M/N}$ (with integral coefficients) which vanishes on $\sigma \cdot x$ for all $\sigma \in M$. Since these vectors span $V$, the linear form vanishes on $V$ and $\chi$ is thus trivial on $\mathrm H_0$. We conclude that $\langle \mathrm e(x)\rangle$ is Zariski dense in $\mathrm H_0$.
\end{proof}

\begin{proof}[Proof of Lemma \ref{l:ZariskiColsureCyclicSL(d,Z)}]
The ring of integers $O_K$ is a free $\mathbb Z$-module of rank $d$. Fixing a basis, one obtains a representation $K^\times \to \GL(d,\mathbb Q)$. This representation comes from a $\mathbb Q$-algebraic morphism $\phi: \mathrm H_0 \to \SL(d)$ which maps $\Lambda$ into $\SL(d,\mathbb Z)$. By Proposition \ref{l:ZariskiDenseCyclicSubgroupTorus}, there exists a cyclic subgroup $\langle y \rangle \in \Lambda$ which is Zariski dense in $\mathrm H_0$. Its image by $\phi$ is a cyclic subgroup of $\SL(d,\mathbb Z)$ which is Zariski dense in a split maximal torus of $\SL(d,\mathbb R)$.
\end{proof}

\subsection{Indiscrete linear hyperbolic groups} 
Let us recall that every simple real rank $1$ Lie group is isogenous to $\mathsf{SO}(n,1)$, $\mathsf{SU}(n,1)$, $\mathsf{Sp}(n,1)$, $n \geq 1$, or to $\textup{F}_4^{(-20)}$ which is the isometry group of the octonionic hyperbolic plane ${\mathbb O}\mathbf{H}^2$. 

We denote by $K_m=\mathsf{Sp}(m)\times \mathsf{Sp}(1)$ the unique up to conjugation maximal compact subgroup of $\mathsf{Sp}(m,1)$. In this subsection, we prove Theorem \ref{t:IndiscretenessIntro}, exhibiting a linear hyperbolic group which does not admit a discrete and faithful representation into any Lie group of rank $1$. 

\indiscreteness*

Note that the subgroup $\langle \Gamma, t\Gamma t^{-1} \rangle$ of the HNN extension $\Gamma \ast_{\langle w \rangle}=\big \langle \Gamma, t \ | \ twt^{-1}=w \big \rangle$ is isomorphic to $\Gamma \ast_{\langle w \rangle} \Gamma$. By Theorem \ref{t:LinearityIntro}, the double $\Gamma \ast_{\langle w \rangle}\Gamma$ admits a faithful representation into $\mathsf{Sp}(k,1)$. Note also that this group is word hyperbolic by the Bestvina--Feighn combination theorem \hbox{(see Theorem \ref{combination}).}

\begin{proof}[Proof of Theorem \ref{t:IndiscretenessIntro}]
Since $k \geq 4$, $\Gamma$ has virtual cohomological dimension at least $16$. Hence, there is no discrete and faithful representation \hbox{$\rho: \Gamma \rightarrow \mathrm{F}_{4}^{(-20)}$}. Indeed, since the symmetric space of $\mathrm{F}_4^{(-20)}$ has dimension 16, such a representation could only exist for $k=4$ and would identify $\Gamma$ to a Zariski dense cocompact lattice in $\mathrm{F}_{4}^{(-20)}$, contradicting Mostow's rigidity \cite{Mostow}. Therefore, since $\mathsf{SO}(m,1) \subset \mathsf{SU}(m,1) \subset \mathsf{Sp}(m,1)$, it is enough to rule out discrete faithful embeddings of the amalgamated product $\Gamma \ast_{\langle w \rangle} \Gamma$ into $\mathsf{Sp}(m,1)$ for every $m \geq 4$. 

Denote by $\Gamma_1$ and $\Gamma_2$ the two copies of $\Gamma$ amalgamated along $\langle w \rangle$ inside $\Gamma \ast_{\langle w \rangle} \Gamma$. Let us also denote by $d_{\mathbb{H}{\bf H}^d}$ the Riemannian metric distance on $\mathbb{H}{\bf H}^d$. Supppose that there exists a discrete and faithful representation \hbox{$\rho:\Gamma \ast_{\langle w \rangle} \Gamma \rightarrow \mathsf{Sp}(m,1)$}. By Theorem \ref{t:SuperRigidityGeodesic}, for $i \in \{1,2\}$, there exists a $\rho|_{\Gamma_i}$-equivariant totally geodesic embedding $f_{i}:\mathbb{H}{\bf H}^k \xhookrightarrow{} \mathbb{H}{\bf H}^m$, and for $i\in \{1,2\}$, $X_{i} \equaldef f_{i}(\mathbb{H}{\bf H}^k)$ is a $\rho(w)$-invariant totally geodesic submanifold of $\mathbb{H}{\bf H}^m$. Therefore, both $X_1$ and $X_2$ contain the axis $L_{\rho(w)}$ of $\rho(w)$. By Fact \ref{fact}, both $X_1$ and $X_2$ contain the unique quaternionic line $L_{\rho(w)}'\subset \mathbb{H}{\bf H}^m$ containing $L_{\rho(w)}$. In particular, $L_{\rho(w)}'$ is contained in $X_1 \cap X_2$. Observe that the action of $\rho(w)$ on $L'_{\rho(w)}$ cannot be cocompact. Hence, we may choose a sequence $(x_n)_{n \in \mathbb{N}}$ in $L'_{\rho(w)}$ such that $\textup{dist}_{\mathbb{H}{\bf H}^m}\big(x_{n}, L_{\rho(w)}\big) \rightarrow \infty$ as $n \rightarrow \infty$. Since $\Gamma_1$ and $\Gamma_2$ are cocompact lattices in $\mathsf{Sp}(k,1)$ and $f_{i}$ is $\rho|_{\Gamma_i}$-equivariant, for every $n \in \mathbb{N}$, we may find $\gamma_n \in \Gamma_1$, $\gamma_{n}' \in \Gamma_2$ and $M>0$ such that $$d_{\mathbb{H}{\bf H}^m}\big(x_{n}, \rho(\gamma_n)K_{m}\big) \leqslant M \ \ \textup{and} \ \ d_{\mathbb{H}{\bf H}^m}(x_{n}, \rho(\gamma_{n}')K_{m}) \leqslant M$$ for every $n \in \mathbb{N}$. In particular, the triangle inequality implies that $$d_{\mathbb{H}{\bf H}^m}\big(\rho(\gamma_{n}^{-1}\gamma_{n}')K_m, K_{m}\big) =d_{\mathbb{H}{\bf H}^m}\big(\rho(\gamma_{n}')K_m, \rho(\gamma_n)K_{m}\big) \leqslant 2M$$ for every $n \in \mathbb{N}$. Since $\rho$ is assumed to have discrete image, we may pass to a subsequence $(k_n)_{n \in \mathbb{N}}$ and assume that $\rho(\gamma_{k_n}^{-1}\gamma_{k_n}')=\rho(\gamma_{k_{n_0}}^{-1}\gamma_{k_{n_0}}')$ or equivalently $\rho(\gamma_{k_{n_0}}\gamma_{k_n}^{-1})=\rho(\gamma_{k_{n_0}}'\gamma_{k_n}'^{-1})$ for every $n \geq n_0$. Since $\rho$ is faithful and $\Gamma_1 \cap \Gamma_2$ is the cyclic group $\langle w \rangle$, there exists $s_n \in \mathbb{N}$ such that $\gamma_{k_{n_0}}\gamma_{k_n}^{-1}=\gamma_{n_0}' \gamma_{k_n}'^{-1}=w^{s_n}$ for $n \geq n_0$. Therefore, there exists $C>0$ with \hbox{$\textup{dist}_{\mathbb{H}{\bf H}^m}(x_{k_n}, L_{\rho(w)})\leq C$} for every $n \in \mathbb{N}$, contradicting the choice of the sequence $(x_n)_{n \in \mathbb{N}}$ above. 

Finally, we conclude that there is no discrete faithful representation $\rho: \Gamma \ast_{\langle w \rangle}\Gamma \rightarrow \mathsf{Sp}(m,1)$. The proof of Theorem \ref{t:IndiscretenessIntro} is complete. \end{proof}

\end{document}